\documentclass[a4paper]{amsart}
\usepackage{graphicx}
\usepackage{amsmath,amsthm,amsfonts,amssymb,mathdots}
\usepackage{afterpage}
\usepackage{latexsym}
\usepackage{subfigure} 
\usepackage[all]{xy}
\xyoption{matrix}

\begin{document}

\renewcommand{\th}{\operatorname{th}\nolimits}
\newcommand{\rej}{\operatorname{rej}\nolimits}
\newcommand{\extto}{\xrightarrow}
\renewcommand{\mod}{\operatorname{mod}\nolimits}
\newcommand{\ul}{\underline}
\newcommand{\Sub}{\operatorname{Sub}\nolimits}
\newcommand{\ind}{\operatorname{ind}\nolimits}
\newcommand{\Fac}{\operatorname{Fac}\nolimits}
\newcommand{\add}{\operatorname{add}\nolimits}
\newcommand{\soc}{\operatorname{soc}\nolimits}
\newcommand{\Hom}{\operatorname{Hom}\nolimits}
\newcommand{\Rad}{\operatorname{Rad}\nolimits}
\newcommand{\RHom}{\operatorname{RHom}\nolimits}
\newcommand{\uHom}{\operatorname{\underline{Hom}}\nolimits}
\newcommand{\End}{\operatorname{End}\nolimits}
\renewcommand{\Im}{\operatorname{Im}\nolimits}
\newcommand{\Ker}{\operatorname{Ker}\nolimits}
\newcommand{\Coker}{\operatorname{Coker}\nolimits}
\newcommand{\Ext}{\operatorname{Ext}\nolimits}
\newcommand{\op}{{\operatorname{op}}}
\newcommand{\Ab}{\operatorname{Ab}\nolimits}
\newcommand{\id}{\operatorname{id}\nolimits}
\newcommand{\pd}{\operatorname{pd}\nolimits}
\newcommand{\A}{\operatorname{\mathcal A}\nolimits}
\newcommand{\C}{\operatorname{\mathcal C}\nolimits}
\newcommand{\D}{\operatorname{\mathcal D}\nolimits}
\newcommand{\X}{\operatorname{\mathcal X}\nolimits}
\newcommand{\Y}{\operatorname{\mathcal Y}\nolimits}
\newcommand{\F}{\operatorname{\mathcal F}\nolimits}
\newcommand{\Z}{\operatorname{\mathbb Z}\nolimits}
\renewcommand{\P}{\operatorname{\mathcal P}\nolimits}
\newcommand{\T}{\mathcal T}
\newcommand{\M}{\mathcal M}
\newcommand{\G}{\Gamma}
\renewcommand{\L}{\Lambda}
\newcommand{\bdot}{\scriptscriptstyle\bullet}
\renewcommand{\r}{\operatorname{\underline{r}}\nolimits}
\newtheorem{lemma}{Lemma}[section]
\newtheorem{prop}[lemma]{Proposition}
\newtheorem{cor}[lemma]{Corollary}
\newtheorem{theorem}[lemma]{Theorem}
\newtheorem{example}[lemma]{Example}
\newtheorem{conjecture}[lemma]{Conjecture}

\theoremstyle{definition}
\newtheorem{remark}[lemma]{Remark}
\newtheorem{definition}[lemma]{Definition}

\title[Ptolemy relations for punctured discs]{Ptolemy relations for
punctured discs}

\date{9 November 2007}

\author[Baur]{Karin Baur}
\address{Department of Mathematics \\
ETH Z\"urich \\
R\"amistrasse 101 \\
8092 Z\"urich \\
Switzerland
}
\email{baur@math.ethz.ch}

\author[Marsh]{Robert J. Marsh}
\address{Department of Pure Mathematics \\
University of Leeds \\
Leeds LS2 9JT \\
England
}
\email{marsh@maths.leeds.ac.uk}


\begin{abstract}
We construct frieze patterns of type $D_N$ with entries which are
numbers of matchings between vertices and triangles of corresponding
triangulations of a punctured disc. For triangulations corresponding to
orientations of the Dynkin diagram of type $D_N$, we show that the numbers
in the pattern can be interpreted as specialisations of cluster variables
in the corresponding Fomin-Zelevinsky cluster algebra.
\end{abstract}

\maketitle

%
\section{Introduction}
\label{sec:introduction}
%

\begin{figure}[htp]
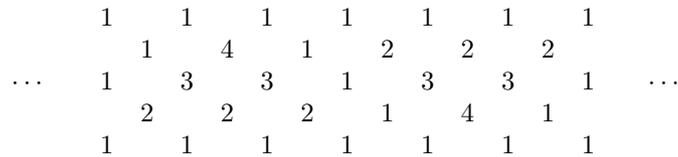

$$
\begin{array}{ccccccccccccccccc}
&& 1 && 1 && 1 && 1 && 1 && 1 && 1 \\
&&& 1 && 4 && 1 && 2 && 2 && 2 & \\
\cdots && 1 && 3 && 3 && 1 && 3 && 3 && 1 && \cdots \\
&&& 2 && 2 && 2 && 1 && 4 && 1 & \\
&& 1 && 1 && 1 && 1 && 1 && 1 && 1
\end{array}
$$
\caption{A Conway-Coxeter frieze pattern of order $6$}
\label{fig:friezeexample}
\end{figure}

Frieze patterns were introduced by Conway and Coxeter
in~\cite{cc73a,cc73b}. Such a pattern consists of a finite number of rows
arranged in an array so that the numbers in the $k$th row sit between
the numbers in the rows on either side. The first and last rows consist
of ones and for every diamond of the form
$\begin{array}{rcl}
 & b \\
a & & d \\
 & c
\end{array}$
the relation $ad-bc=1$ must be satisfied. The order $N$ of the pattern is one
Coxeter and Conway associated a frieze pattern of order $N$
(i.e. with $N-1$ rows) to each triangulation of a regular polygon with
$N$ sides and showed that every frieze pattern arises in this way. For
an example see Figure~\ref{fig:friezeexample}.

In~\cite{calderochapoton06} and~\cite{propp}, the authors consider
frieze patterns arising from Fomin-Zelevinsky cluster algebras of
type $A_N$~\cite{fominzelevinsky02,fominzelevinsky03}. Motivated by the
description~\cite{fst06} of cluster algebras of type $D_N$ in terms of
(tagged) triangulations of a disc $S_{\odot}$ with a single puncture and
$N$ marked points on the boundary we associate a frieze pattern to every
such triangulation. Note that the type $D_N$ case
of~\cite{fst06} has been described in detail by Schiffler~\cite{schiffler06}
in his study of the corresponding cluster category.

Each number in our frieze pattern is the cardinality of a
set of matchings of a certain kind between vertices of the triangulation
and triangles in it, so our result can be regarded as a generalisation
of a result of~\cite{carrollprice} (see~\cite{propp}) giving the
numbers in a Conway-Coxeter frieze pattern in terms of numbers of
perfect matchings for graphs associated to the corresponding triangulation
of an unpunctured disc.

We show further that, in the cases where the triangulation has a particularly
nice form (corresponding to an orientation of the Dynkin diagram of
type $D_N$), the numbers in our frieze pattern can be interpreted as
specialisations of cluster variables of the cluster algebra of type $D_N$,
generalising similar results in type $A_N$ obtained by Propp~\cite{propp}
(for arbitrary triangulations).

We remark that in independent work Musiker~\cite{musiker} has recently
associated a frieze pattern to the initial bipartite seed of a
cluster algebra of types $A,B,C$ or $D$ in which the entries are
numbers of perfect matchings of associated graphs given by tilings
(in fact via weightings of edges, Musiker is also able to describe the
numerators of the corresponding cluster variables).
In type $D_N$ this corresponds to one triangulation in our picture
(up to symmetry).

\begin{figure}
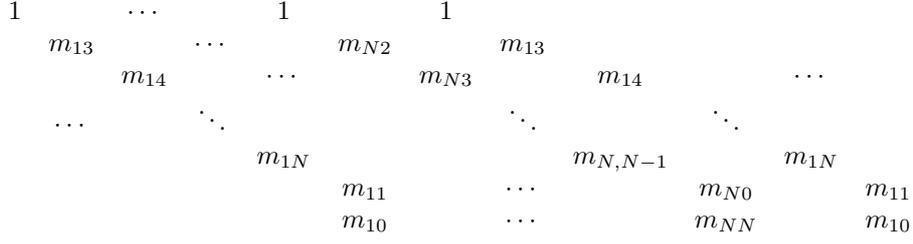

$$\begin{array}{ccccccccccccc}
1 && \cdots && 1 && 1 &&&&&& \\
&m_{13}&      &\cdots&         &m_{N2} &         &m_{13}    &        &    &        &       &        \\
&      &m_{14}&      &\cdots   &         &m_{N3} &          &m_{14} &    &       \cdots  &  &        \\
&\cdots&         &\ddots   &         &        &          &\ddots  &     &\ddots       &       &        \\
&      &      &         &m_{1N}  &        &          &        &m_{N,N-1}&             &m_{1N}&        \\
&      &      &         &         &m_{11}&          &\cdots  &         &m_{N0}    &       &m_{11}\\
&       &      &         &         &m_{10}&          &\cdots  &         &m_{NN}    &       &m_{10}
\end{array}
$$
\caption{A frieze pattern of type $D_N$, $N$ even.}
\label{fig:dpatterneven}
\end{figure}

\begin{figure}
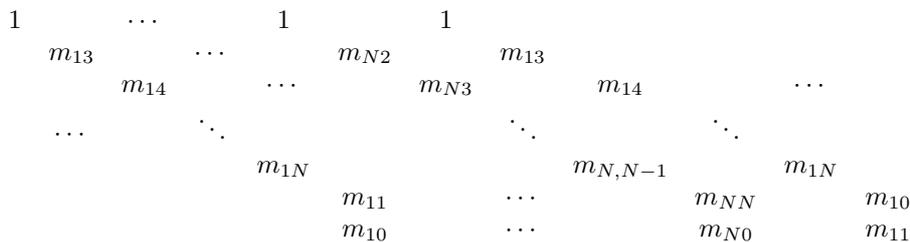

$$\begin{array}{ccccccccccccc}
1 && \cdots && 1 && 1 &&&&&& \\
&m_{13}&      &\cdots&         &m_{N2} &         &m_{13}    &        &    &        &       &        \\
&      &m_{14}&      &\cdots   &         &m_{N3} &          &m_{14} &    &       \cdots  &  &        \\
&\cdots&         &\ddots   &         &        &          &\ddots  &     &\ddots       &       &        \\
&      &      &         &m_{1N}  &        &          &        &m_{N,N-1}&             &m_{1N}&        \\
&      &      &         &         &m_{11}&          &\cdots  &         &m_{NN}    &       &m_{10}\\
&       &      &         &         &m_{10}&          &\cdots  &         &m_{N0}    &       &m_{11}
\end{array}
$$
\caption{A frieze pattern of type $D_N$, $N$ odd.}
\label{fig:dpatternodd}
\end{figure}

\begin{figure}
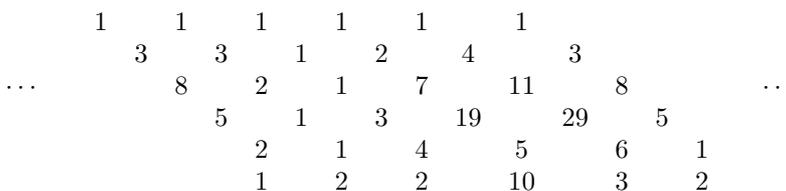

$$\begin{array}{ccccccccccccccccccc}
&&1 && 1 && 1 && 1 && 1 && 1 &&&&&& \\
&&& 3 && 3 && 1 && 2 && 4 && 3 &&&&& \\
\cdots &&&& 8 && 2 && 1 && 7 && 11 && 8 && && \cdots \\
&&&&& 5 && 1 && 3 && 19 && 29 && 5 &&& \\
&&&&&& 2 && 1 && 4 && 5 && 6 && 1 &&\\
&&&&&& 1 && 2 && 2 && 10 && 3 && 2 &&
\end{array}
$$
\caption{A frieze pattern of type $D_5$.}
\label{fig:D5examplefrieze}
\end{figure}

We now go into more detail.
Arcs in triangulations of $S_{\odot}$ can be indexed by ordered pairs $i,j$ of
(possibly equal) boundary vertices $i,j$ together with
unordered pairs $i,0$ consisting of a vertex on the boundary together with
the puncture, with the proviso that the pairs $i,i+1$ (where $i+1$ is
interpreted as $1$ if $i=N$) are not allowed. We denote the arcs by
$D_{ij}$ and $D_{i0}$ respectively. Let $m_{ij}$ and $m_{i0}$ be
positive integers associated to such arcs.
We define a \emph{frieze pattern of type $D_N$} to
be an array of numbers of the form in Figure~\ref{fig:dpatterneven} if
$N$ is even or of the form in Figure~\ref{fig:dpatternodd} if $N$ is odd.
Note that in the odd case, every other occurrence of the fundamental pattern
has the bottom two rows flipped.
The following relations must be satisfied for all boundary vertices
$i,j$:
\begin{eqnarray}
m_{ij}\cdot m_{i+1,j+1} & = & m_{i+1,j}\cdot m_{i,j+1} + 1, 
\mbox{ provided $j\neq i$ or $i+1$}; \label{eqn:relation1} \\
m_{i,i-1}\cdot m_{i+1,i} & = & m_{i+1,i-1}\cdot m_{ii}\cdot m_{i0}+1;
\label{eqn:relation2} \\
m_{ii}\cdot m_{i+1,0} & = & m_{i+1,i} + 1; \label{eqn:relation3} \\
m_{i0}\cdot m_{i+1,i+1} & = & m_{i+1,i} + 1. \label{eqn:relation4}
\end{eqnarray}
We refer to these relations as the \emph{frieze relations}. For an example
of a frieze pattern of type $D_5$ see Figure~\ref{fig:D5examplefrieze}.

In Definition~\ref{def:matchingnumbers} we will explain how to associate
numbers of matchings $m_{ij}$ and $m_{i0}$ to the arcs of $S_{\odot}$ (with
$N$ marked points on the boundary).

\noindent {\bf Theorem A.}
\emph{If the matching numbers $m_{ij}$ are arranged as above,
then they form a frieze pattern of type $D_N$.}

Let $\mathcal{A}$ be a cluster algebra~\cite{fominzelevinsky02} of type
$D_N$. We consider the case in which all coefficients
are set to $1$. Let $x_{ij}$ be the cluster variable corresponding
to the arc in $S_{\odot}$ with end-points $i$ and
$j$~\cite{fst06,schiffler06}. 
Fix a seed $(\mathbf{x},B)$ of $\mathcal{A}$, so that $\mathbf{x}$ is
a cluster and $B$ is a skew-symmetric matrix. Let $Q$ be the
quiver associated to $B$. Let $u_{ij}$ denote the integer obtained
from $x_{ij}$ when the elements of $\mathbf{x}$ are all specialised to $1$.

\noindent \textbf{Theorem B.}
\emph{Suppose that the quiver $Q$ of the seed $(\mathbf{x},B)$ is
an orientation of the Dynkin diagram of type $D_N$.
Let $D_{ij}$ be an arc in $S_{\odot}$. Then $m_{ij}=u_{ij}$.}

In Section 2 we introduce the necessary notation and definitions,
describing the matching numbers referred to above in detail. In
Section 3 we recall the properties of the numbers in Conway-Coxeter
frieze patterns that we need. In Section 4 we prove Theorem A
and in Section 5 we prove Theorem B.

%
\section{Notations and definitions}
\label{sec:notation}
%

We consider ``triangulations''
of bounded discs with a number
of marked points on the boundary and up to
one marked point in the interior (a puncture).
Such a disc is a special case of a bordered surface 
with marked points as in
recent work~\cite{fst06} of Fomin, Shapiro and Thurston
- they consider surfaces with an arbitrary number of boundary
components and punctures.
We use the following notation:
$S$ denotes a bounded disc with $N$ marked points
(or boundary vertices) on the boundary, $N\ge 3$.
We add a subscript $_{\odot}$
if the disc has a puncture.
We will usually label the points on the boundary clockwise around the
boundary with $1,2,\dots,N$ and denote the puncture by $0$.

For any such disc
let $B_{ij}$ denote the boundary arc
from the vertex $i$ up to the vertex $j$, {\em going clockwise}.
So if $j=i+1$, $B_{i,i+1}$ denotes the boundary component
from $i$ to $i+1$ including the two marked points.
({\em Marked points are always taken $\mod N$}).
If the disc has no puncture, we assume $i\neq j$.
If it has a puncture then $i=j$ is allowed:
$B_{ii}$ denotes the whole boundary, starting and ending at $i$.

We consider arcs {\em inside a disc} up to isotopy:
\begin{definition} \label{def:arc}
An {\em arc} $D$ of a disc is an element
of the isotopy class of curves whose
endpoints are marked points of the disc. The arc is not
allowed to intersect itself.
Furthermore we require that the interior of the
arc is disjoint from the boundary of the disc and that it does
not cut out an unpunctured monogon or digon.
\end{definition}

In other words, an arc $D$ is of the form
$D_{ij}$ with endpoints $i,j$ where
for an unpunctured disc $1\le i,j\le N$,
$j\notin\{i-1,i,i+1\}$, and
for a punctured disc
$i,j\ \in\{1,\dots,N\}\cup\{0\}$, $j\neq i+1$ and
$i,j$ not both equal to $0$.
So arcs of the form $D_{i0}$, $D_{i,i-1}$,
$D_{ii}$ only occur in punctured discs. The latter is an
arc with common endpoints, winding around the puncture once.
Note that for $1\le i\neq j\le N$ we have
$D_{ij}=D_{ji}$ if and
only if the disc has no puncture. If the disc has a puncture
$0$, and if $i\neq j$, $D_{ij}$ winds partly around $0$, going
clockwise
from $i$ to $j$ and $D_{ji}$ winds around the puncture going
clockwise
from $j$ to $i$, i.e. anti-clockwise
from $i$ to $j$.
Note also that by definition, a boundary arc $B_{i,i+1}$ is not
considered to
be an arc.

Using the arcs, the discs can be triangulated: 

\begin{definition}
A {\em triangulation} $\T$ of a disc $S$ (or $S_{\odot}$)
is a maximal collection $\T$ of pairwise non-intersecting arcs
of $S$ (of $S_{\odot}$).
\end{definition}

\begin{figure}
\begin{center}
\includegraphics[scale=.4]{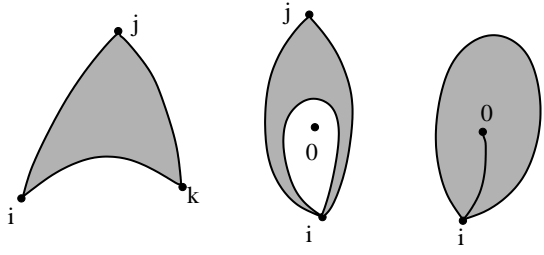}
\end{center}
\caption{The three types of triangles}
\label{fig:triangles}
\end{figure}

The arcs of the triangulation divide the disc into a collection of
disjoint triangles. We will often call them the {\em triangles} of
$\T$.
There are three types of triangles appearing, as shown
in Figure~\ref{fig:triangles}.
If the disc has no puncture, the only triangles appearing have
three vertices and three arcs as sides.
If the disc has a puncture, there are two more types of triangles,
one with sides $D_{ij}$, $D_{ji}$ and $D_{ii}$ ($i,j\neq 0$)
and one with
the corresponding sides obtained by setting $j=0$.
The latter is called a {\em self-folded} triangle.
These three types correspond to the triangles of label A$_1$,
A$_2$ and A$_4$ of Burman, cf.~\cite{burman99} and are called
{\em ideal} triangles in~\cite{fst06} - our notion of triangulation
is thus an {\em ideal triangulation} in~\cite{fst06}.

The number of arcs in a triangulation is an invariant of the disc.
It is clear that any triangulation of $S$ has $N-3$ arcs and $N-2$
triangles
whereas triangulations of punctured discs $S_{\odot}$ have $N$ arcs
and $N$ triangles.

We often need to refer to a subset $V_{kl}$ of the vertices of the
disc:
\begin{definition}
For $1\le k,l\le N$ we write $V_{kl}:=V(B_{k,l})$
to denote the vertices of the
boundary arc $B_{k,l}$. In other words, we are referring to the
vertices
from $k$ to $l$, going clockwise. In case $k>l$, the set $V_{kl}$
consists of the
vertices $k,k+1,\dots,N,1,2,\dots,l$.
\end{definition}

We will also need to consider truncations
$(S_{ij})_{\odot}$ and $S_{ij}$ of the punctured disc
$S_{\odot}$, defined as follows:
\begin{definition}\label{def:truncation}
Let $j\neq i+1$.
Together with $D_{ij}$, the boundary arc
$B_{ji}$ forms a
(smaller) punctured disc which we denote by
$(S_{ij})_{\odot}$. Its vertices
are $V_{j,i}$. In particular, $j$ is the clockwise
neighbour of $i$ in $(S_{ij})_{\odot}$.

On the other hand, $D_{ij}$ and $B_{ij}$ form
an unpunctured disc which we denote by $S_{ij}$.
Its vertices are $V_{i,j}$, where $i$ is the
clockwise neighbour of $j$ in $S_{ij}$.
\end{definition}

In $S_{ij}$ and $(S_{ij})_{\odot}$,
some of the triangles of
$\T$ may be cut open into several different regions.
If $D_{ij}$ is not an arc of the
triangulation, then there are triangles of $\T$
which are crossed by $D_{ij}$.
If $S_{ij}$ contains two connected components of
such a triangle we say that $D_{ij}$ {\em splits} the triangle.

\begin{remark}
If $D_{ij}$ is an arc of the triangulation, it does not cut across
any triangles. Hence no split triangles appear and the arcs
of $\T$ induce triangulations of $(S_{ij})_{\odot}$ and of $S_{ij}$.
\end{remark}

\begin{definition}
Let $S_{\odot}$ be a punctured disc with triangulation $\T$. \\
We denote by $\T|_{S_{ij}}$ the subdivision of $S_{ij}$ into
triangles and regions obtained from restricting the
triangulation $\T$ to $S_{ij}$. \\
In the same way, $\T|_{(S_{ij})_{\odot}}$ denotes the subdivision of
$(S_{ij})_{\odot}$ into triangles and regions obtained from
restricting $\T$ to $(S_{ij})_{\odot}$.
\end{definition}

\begin{figure}
\begin{center}
\subfigure[no splitting]
{\includegraphics[scale=.4]{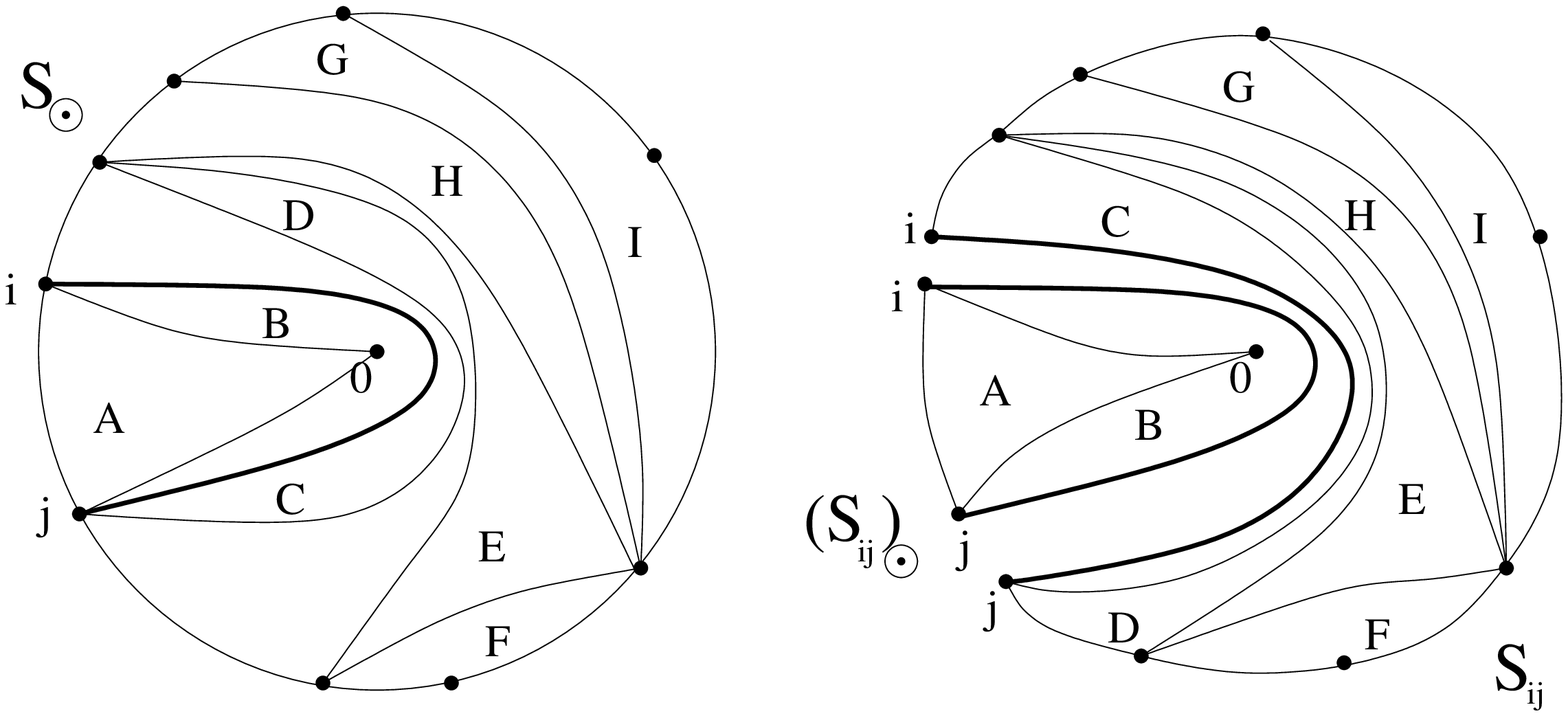}\label{fig:Non-split}}
\hfill
\subfigure[splitting occurs]
{\includegraphics[scale=.4]{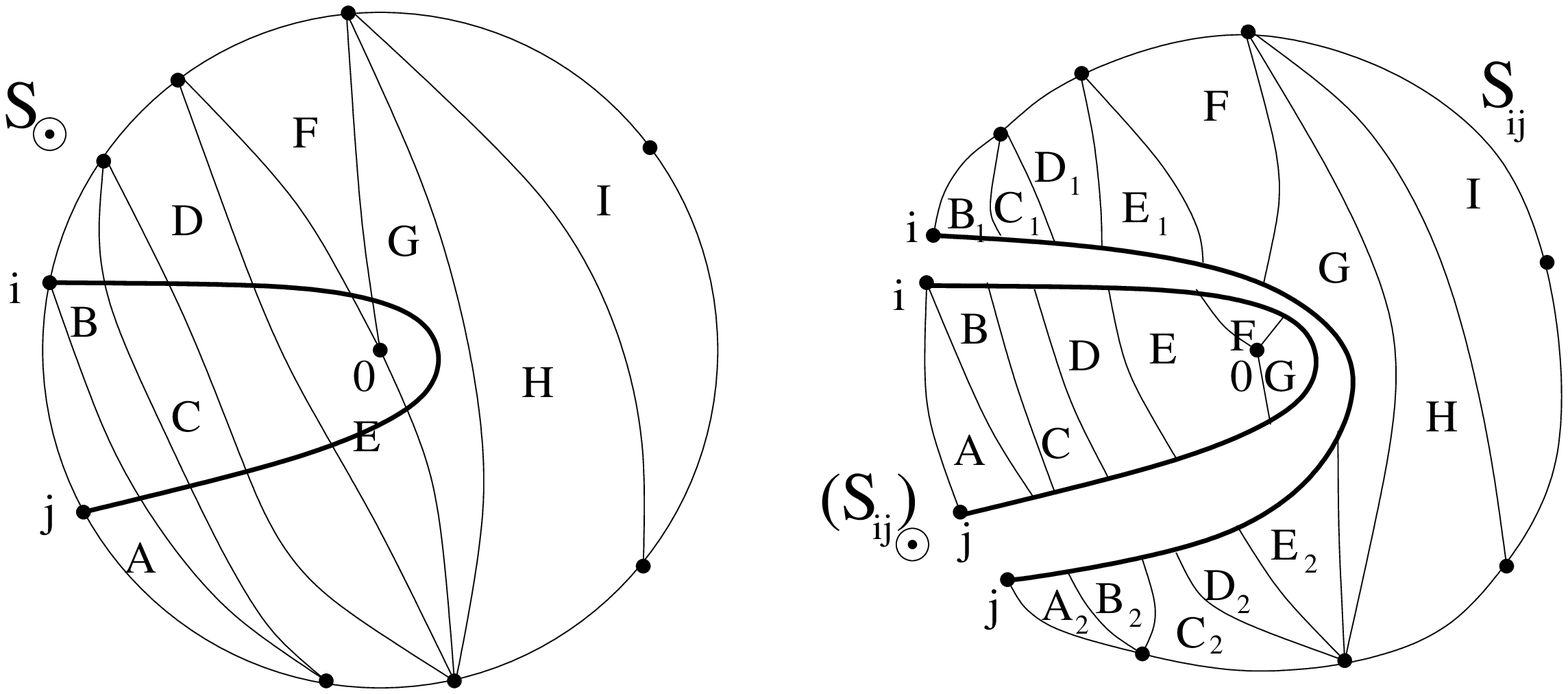}\label{fig:split}}
\end{center}
\caption{Restricted triangulations through truncation by the arc
$D_{ij}$}
\label{fig:truncating}
\end{figure}

Figure~\ref{fig:truncating} presents two examples, one where
no triangles are split and one where $D_{ij}$ splits triangles.
Note that by the observations above, $\T|_{S_{ij}}$ and
$\T|_{(S_{ij})_{\odot}}$
are triangulations of the corresponding discs if and only if
$D_{ij}\in\T$.

\begin{remark}\label{rem:splitting}
Consider the arc $D_{ij}$ (where $j\neq i+1$) in a
punctured disc $S_{\odot}$ with triangulation $\T$. If there is
a $k$ in $V_{j,i}$ such that $\T$ contains the central
arc $D_{k0}$, then $\T|_{S_{ij}}$ does not contain any pair of split
triangles.
\end{remark}

Given a triangulation of any disc (punctured or unpunctured) we
are interested in the ways to allocate triangles to a subset $I$
of the vertices. That is we are interested in matchings between a
subset of the vertices and a subset of the triangles. In the case
of a truncated version of a disc, we want to allocate triangles
{\em and} regions of the truncated triangulation to vertices.

\begin{definition}\label{def:matching}
Let $S$ and $S_{\odot}$ be discs with vertices
$\{1,2,\dots,N\}$, $S_{\odot}$ with puncture $0$. Fix a 
triangulation $\T$ of the disc, let $I$ be a subset of the vertices
$\{1,\dots,N\}$. 

(i) A {\em matching (for $S$) between $I$ and $\T$}
is a way to allocate to each vertex $i\in I$ a triangle of $\T$
incident with $i$ in such a way that no triangle of $\T$ is
allocated to more than one vertex.

(ii) A {\em matching (for $S_{\odot}$) between $I\,\cup\,\{0\}$
and $\T$} is a way to allocate $k+1$ triangles
of $\T$ to the vertices of
$I$ and to the puncture in the same way as above.

(iii) Given an arc $D_{ij}$ in $S_{\odot}$ let $S_{ij}$ be the
truncated disc as in Definition~\ref{def:truncation} and let
$I\subset V_{i,j}$ be a subset of the vertices between $i$
and $j$. Then a {\em matching (for $S_{ij}$) between
$I$ and $\T|_{S_{ij}}$} is a way to allocate $k=|I|$
triangles or regions of $\T|_{S_{ij}}$ to the vertices of $I$
as above.
\end{definition}

\begin{remark}
Matchings as in Definition~\ref{def:matching} can regarded as
\emph{complete matchings} for the corresponding bipartite graphs
(cf.~\cite[11-1]{grimaldi85}).
\end{remark}

Now if $\T$ is a triangulation of a disc (with or without puncture),
we want to refer to the set of all matchings between a subset of
the vertices of the disc and the triangles of $\T$:
\begin{definition}\label{def:matching-set}
Let $\T$ be a triangulation of a disc and $I$ be a subset of
the vertices $\{1,\dots,N\}$ or of  $\{1,\dots,N\}\cup\{0\}$ in case
the disc has a
puncture. Then we set
$\M(I,\T)$ to be the set of all matchings between $I$ and $\T$.
\end{definition}

\begin{figure}
\subfigure
{\includegraphics[scale=.5]{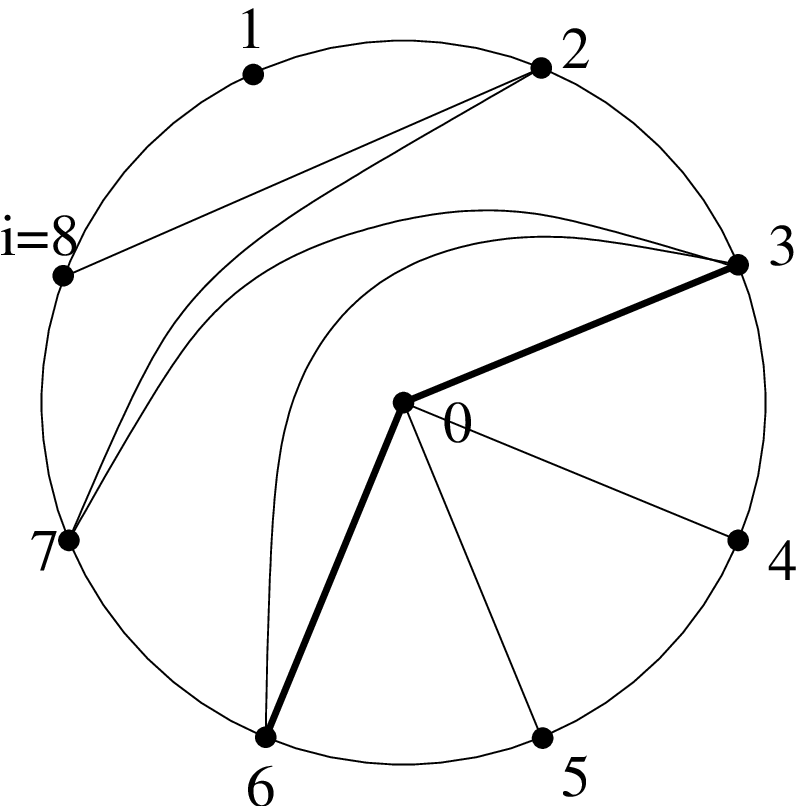}\label{fig:P_i_2central}}
\hspace{1cm}
\subfigure
{\includegraphics[scale=.5]{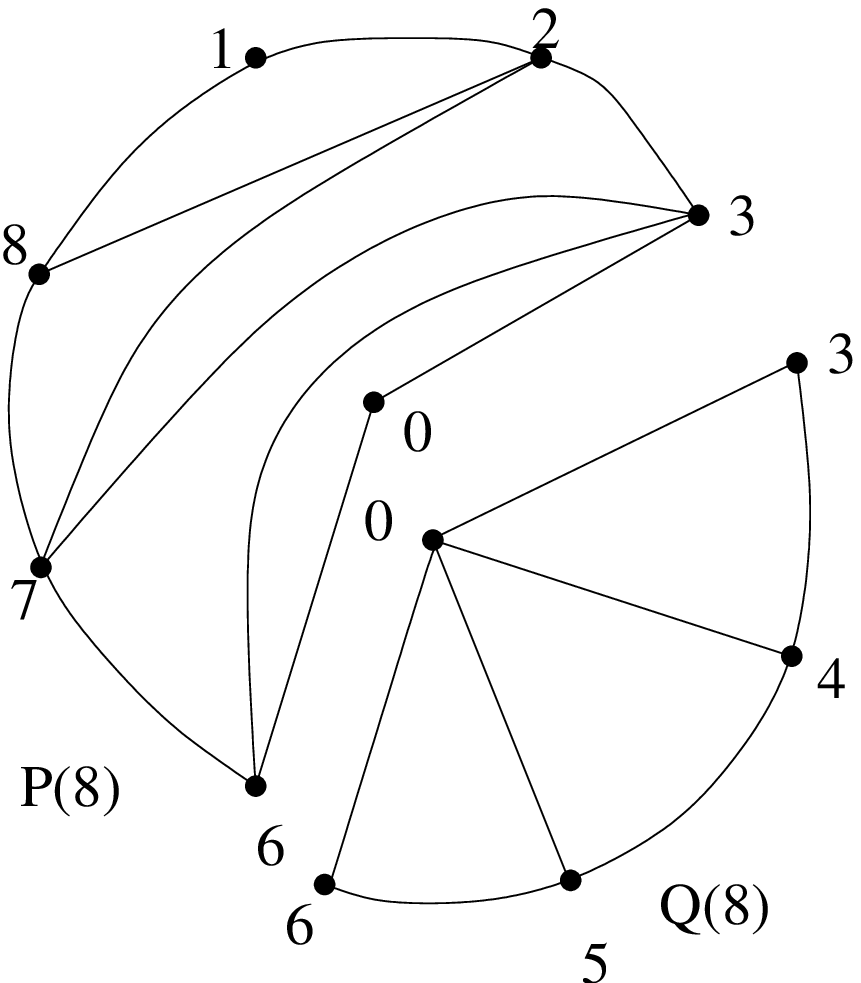}\label{fig:P_i_split}}
\caption{Type (i), $i=8$}\label{fig:(i)-split}
\end{figure}

\begin{figure}
\subfigure
{\includegraphics[scale=.5]{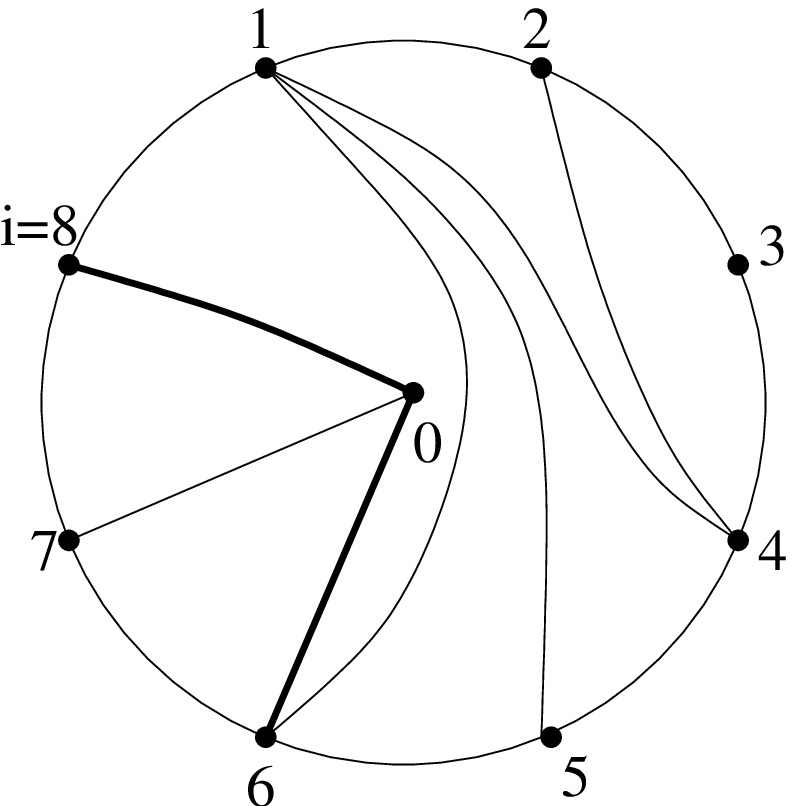}\label{fig:i0-in-T2}}
\hspace{1cm}
\subfigure
{\includegraphics[scale=.5]{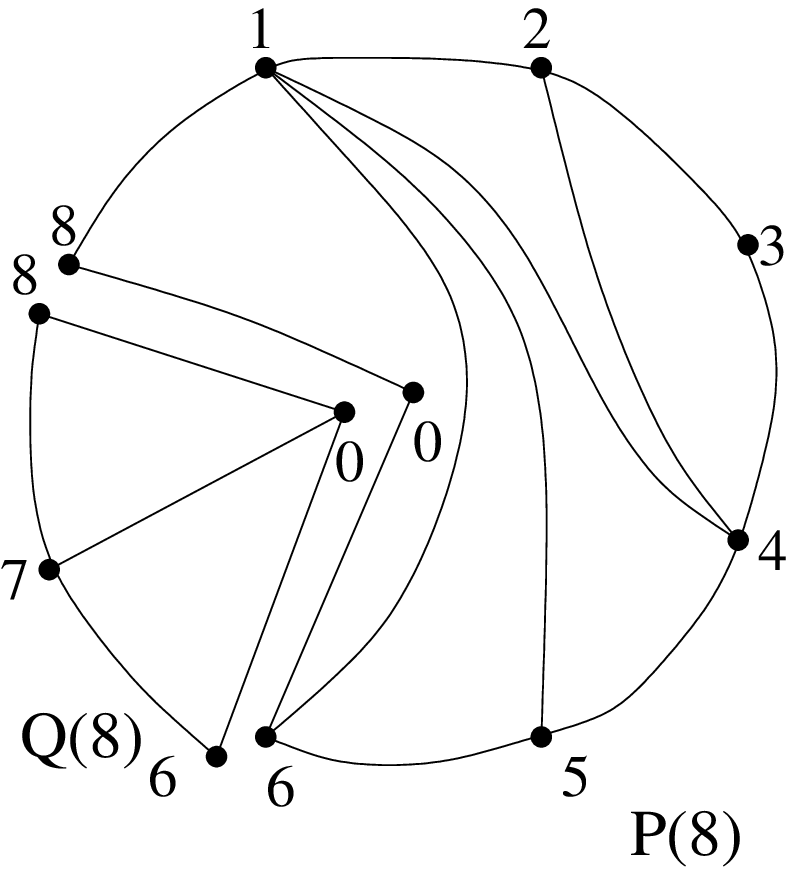}\label{fig:i0k0-split}}
\caption{Type (ii), $i=8$}\label{fig:(ii)-split}
\end{figure}

\begin{figure}
\subfigure
{\includegraphics[scale=.5]{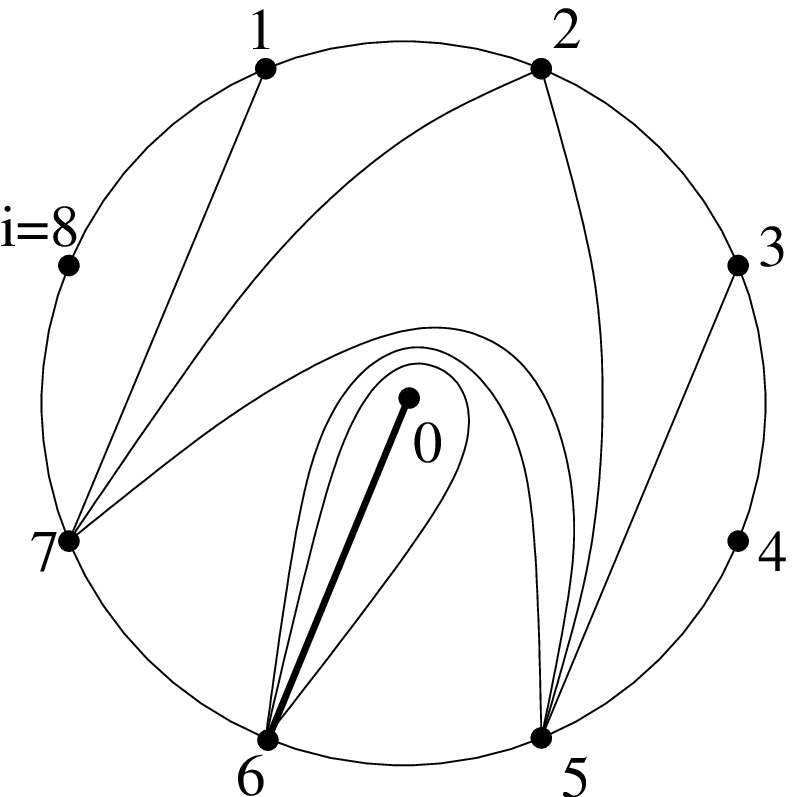}\label{fig:P_i_1central}}
\hspace{1cm}
\subfigure
{\includegraphics[scale=.5]{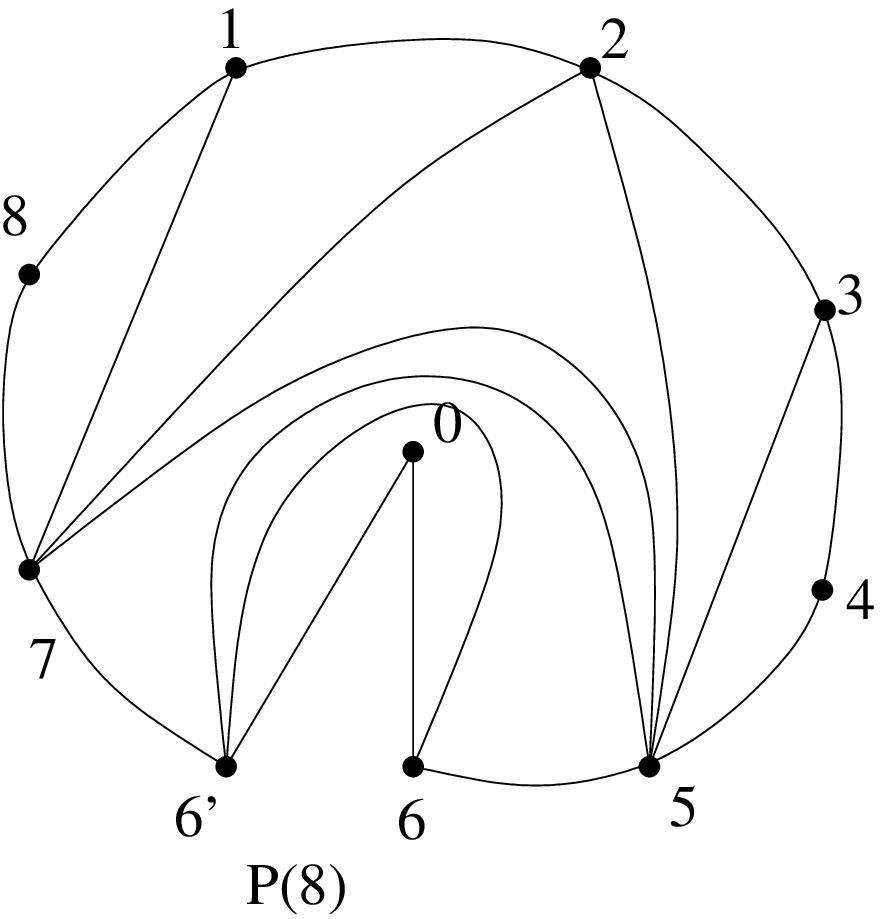}\label{fig:P_i_unfold}}
\caption{Type (iii), $i=8$}\label{fig:(iii)-split}
\end{figure}

\begin{figure}
\subfigure
{\includegraphics[scale=.5]{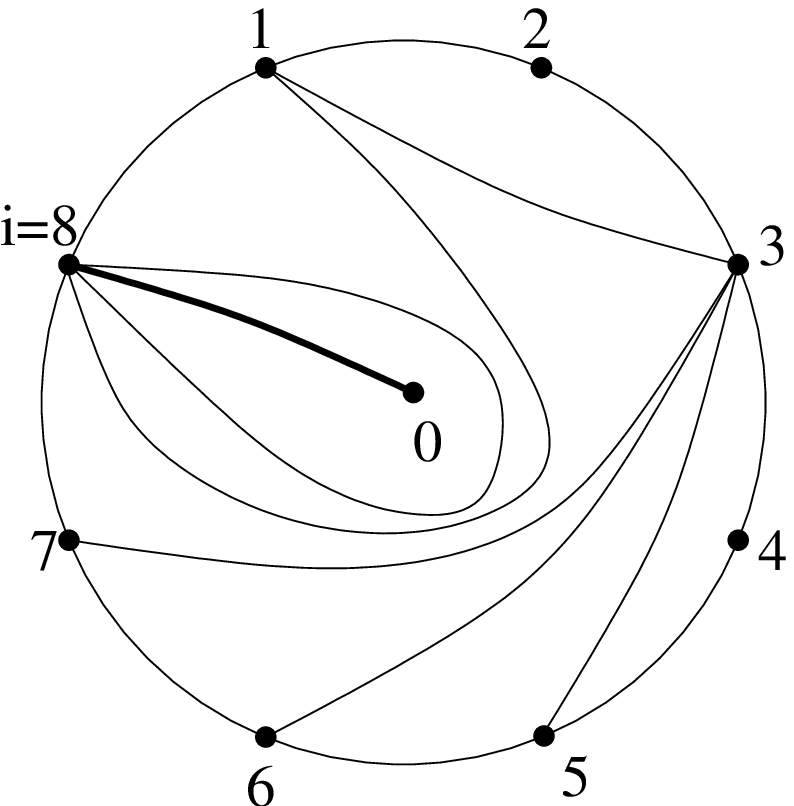}\label{fig:i0-in-T1}}
\hspace{1cm}
\subfigure
{\includegraphics[scale=.5]{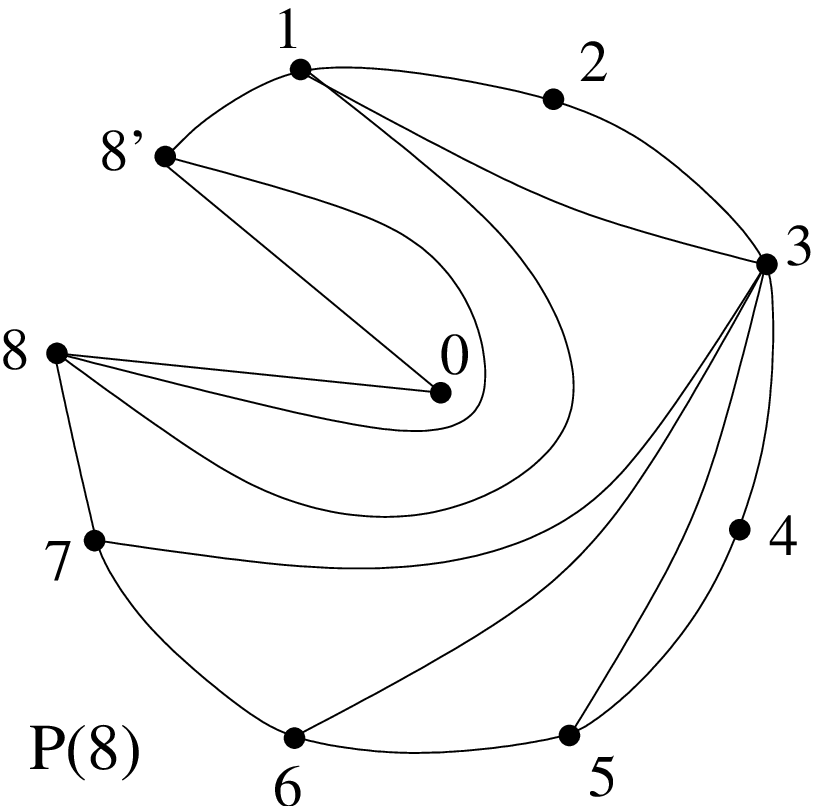}\label{fig:i0-unfold}}
\caption{Type (iv), $i=8$}\label{fig:(iv)-split}
\end{figure}

\begin{remark}
Let $S_{\odot}$ be a punctured disc with $N$ marked points
on the boundary and a fixed triangulation $\T$, and 
let $i\in\{1,2,\dots,N\}$. Then we distinguish
four cases. They are illustrated on the left hand
sides of Figures~\ref{fig:(i)-split},
~\ref{fig:(ii)-split},~\ref{fig:(iii)-split} and
~\ref{fig:(iv)-split}:

(i) $\T$ contains at least two central arcs
and $D_{i0}\notin\T$;

(ii) $\T$ contains at least two central arcs and
$D_{i0}$ $\in\T$;

(iii) There is exactly one central arc $D_{k0}$ in the
triangulation, $i\neq k$;

(iv) The only central arc of $\T$ is $D_{i0}$.
\end{remark}

We want to associate certain unpunctured discs
$P(i)$ and $Q(i)$ to the vertex $i$, in a similar way
as we use an arc $D_{ij}$ to define the truncated discs
$S_{ij}$ and $(S_{ij})_{\odot}$ obtained from $S_{\odot}$.
Here we use one or two central arcs to cut open $S_{\odot}$.

\begin{definition} \label{def:piqi}
(i)
Let $i$ be a marked point of $S_{\odot}$ and $j,k$ 
boundary vertices 
with $i\in V_{kj}$ such that $D_{j0},D_{k0}\in\T$
and such that there is no $l\in V_{kj}$ other than $k$, $j$ 
with $D_{l0}\in\T$.
Then we let $P(i)$ be the unpunctured disc with boundary
$B_{kj}\cup D_{j0}\cup D_{k0}$. 
On the other hand we let $Q(i)$ be the unpunctured disc with
boundary
$B_{jk}\cup\ D_{k0}\cup\ D_{j0}$. Figure~\ref{fig:(i)-split}
illustrates this.

(ii) Let $i$ be  a boundary vertex with $D_{i0}\in \T$ and 
let $j$ be the nearest clockwise neighbour of $i$ with $D_{j0}\in\T$. 
Then we let $P(i)$ be the unpunctured disk with boundary 
$B_{ij}\cup D_{j0}\cup D_{i0}$ and $Q(i)$ be the unpunctured 
disc with boundary  $B_{ji}\cup\ D_{i0}\cup\ D_{j0}$. This 
is illustrated in Figure~\ref{fig:(ii)-split}.

(iii) If the only central arc of
$\T$ is $D_{k0}$ then we define $P(i)$ to be the unpunctured disc obtained 
by cutting up the disc at $k$, with boundary $B_{kk}\cup D_{k0}\cup D_{k0}$. 
We obtain an
additional marked point on the boundary. We denote the anti-clockwise
neighbour of $0$ by $k$ and the clockwise neighbour of $0$ by $k'$.
An example is presented in Figure~\ref{fig:(iii)-split}.

(iv) If the only central arc of
$\T$ is $D_{i0}$ then we define $P(i)$ to be the unpunctured disc obtained 
by cutting up the disc at $i$, with boundary $B_{ii}\cup D_{i0}\cup D_{i0}$. 
As before, we obtain an additional marked point on the boundary. 
We denote the anti-clockwise
neighbour of $0$ by $i$ and the clockwise neighbour of $0$ by $i'$.
An example is presented in Figure~\ref{fig:(iv)-split}.
\end{definition}

In other words, in the first case we cut the disc $S_{\odot}$ along
two
central arcs and $Q(i)$ is the complement to $P(i)$.
In the latter case we unfold the disc along the only central
arc. In particular, $Q(i)$ is not defined in this case.

\begin{remark}
The vertex $i$ and the puncture $0$ are neighbours
in $P(i)$ if
and only if the arc $D_{i0}$ is in the triangulation.
\end{remark}

\begin{definition}
For $P(i)$ and $Q(i)$ as defined above, we let $\T|_{P(i)}$
and $\T|_{Q(i)}$ be the triangulations of $P(i)$ and $Q(i)$
obtained from $\T$ by only considering the
triangles in $P(i)$ and in $Q(i)$ respectively.
Then $\T|_{P(i)}$ and $\T|_{Q(i)}$
are triangulations of unpunctured discs.
\end{definition}

With the notation introduced in Definition~\ref{def:matching-set}
we can now associate numbers $m_{ij}$ to a disc with a fixed
triangulation.

\begin{definition} \label{def:matchingnumbers}
Let $S_{\odot}$ be a punctured disc with vertices
$\{1,2,\dots,N\}$ and puncture $0$.
Fix a triangulation $\T$ of $S_{\odot}$. Let $i,j \in\{1,\dots,N\}$,
$j\neq i,i+1$.
\begin{itemize}
\item[(i)] 
Let $\M_{ij}:=\M(V_{i+1,j-1},\T|_{S_{ij}})$, i.e. 
the set of matchings between $V_{i+1,j-1}$ and $\T|_{S_{ij}}$.
\item[(ii)] 
Let $\M_{i0}:=\M(V(P(i)\setminus\{i,0\}),\T\mid_{P(i)})$.
\item[(iii)] 
Let $\M_{ii}:=\M(V_{i+1,i-1}\cup\{0\},\T)$.
\end{itemize}
We set $m_{i,i+1}=1$, $m_{ij}=|\M_{ij}|$,
$m_{i0}=|\M_{i0}|$ and $m_{ii}=|\M_{ii}|$.
\end{definition}

%
\section{The unpunctured case}
\label{sec:unpuncturedcase}
%

In this section we recall the frieze patterns of Conway and Coxeter
and describe various interpretations of their entries and properties
that they satisfy. We will phrase these results in terms of our set-up
above.

In~\cite{cc73a},~\cite{cc73b}, Conway and Coxeter studied frieze patterns of
positive numbers. These are patterns of positive integers arranged
in a finite number of rows where the top and bottom rows consist of $1$'s.
The entries in the second row appear between the entries in the first
row, the entries in the third row appear between the entries in the second
row, and so on. Also, for every {\em diamond} of the form
$\begin{array}{rcl}
 & b \\
a & & d \\
 & c
\end{array}$
the relation $ad-bc=1$ must be satisfied. The order $N$ of the pattern is one
more than the number of rows. For an example see
Figure~\ref{fig:friezeexample}.

Conway and Coxeter proved in~\cite{cc73a,cc73b} that the second row
of a frieze pattern of order $n$ is equal to the sequence of numbers of
triangles at the vertices of a triangulation of a disc with $N$ marked
points.
Start with a vertex $i$ of $S$ and label it $0$. Whenever $i$ is connected
to another vertex $j$ by an edge (including boundary edges) of $\T$,
label $j$ by $1$. Furthermore, if $\triangle$ is a triangle with exactly
two vertices which have been labelled already, label the third vertex with
the sum of the other two labels. Iterating this procedure, we obtain
labels for all vertices of $S$. The labels clearly depend on $i$
and the label obtained at $j$ (for $j\neq i$) is denoted by $(i,j)$.

\begin{remark} \label{rem:unpuncturedone}
Note that it is clear from the definition that $(i,i+1)=1$ for all $i$
and that $(i,j)=1$ whenever $i$ and $j$ are the two ends of an arc of $\T$.
\end{remark}

The frieze pattern corresponding to $\T$ can then be described by
displaying the numbers $(i,j)$ in the following way:

$$
\begin{array}{ccccccccc}
(1,2) && (2,3) && \cdots && \cdots && (N-1,N) \\
& (1,3) && \cdots && \cdots && (N-2,N) & \\
&& \ddots &&  && \iddots && \\
&&& (1,N-1) && (2,N) &&& \\
&&&& (1,N) &&&&
\end{array}
$$

This fundamental region is then repeated upside down to the right
and left, then the right way up on each side, and so on.
Broline, Crowe and Isaacs have given the following interpretation of the
numbers $(i,j)$.

\begin{theorem} \cite[Theorem 1]{bci74} \label{thm:broline}
Let $S$ be a disc with $N$ marked points on the boundary
with triangulation $\T$.
Let $i,j$ be distinct marked points.
Then: \\
(a) We have that $(i,j)$ is equal to $|\M(V_{i+1,j-1},\T)|$, i.e.\ 
the number of matchings between $\T$ and $i+1,\ldots ,j-1$.\\
(b) We have that $(i,j)=|\M(V_{j+1,i-1},\T)|$. \\
In particular, $(i,j)=(j,i)$.
\end{theorem}

\begin{definition}
Let $i,j$ be a distinct pair of marked points on the boundary of $S$,
and let $\T$ be a triangulation of $S$.
Let $n_{ij}=M(\{1,\dots,N\}\setminus\{i,j\},\T)$, i.e.\ 
the number of matchings between
$\{1,2,\dots,N\}\setminus\{i,j\}$ and $\T$.
\end{definition}

Carroll and Price have shown that the numbers
$(i,j)$ coincide with the $n_{ij}$ defined above. This is discussed
in~\cite{propp}.

\begin{theorem} \cite{carrollprice} \label{thm:bracketequalsn}
With the notation above, $(i,j)=n_{ij}$ for any pair of distinct
vertices $i,j$.
\end{theorem}

We shall freely use this result in the sequel, in particular noting
that the $n_{ij}$ have the properties of the $(i,j)$ that we have
described above.

Propp reports that one of the main steps in the proof of Carroll and
Price is the following, which is a direct consequence of
the \emph{Condensation Lemma} of Kuo~\cite[Theorem 2.5]{kuo04}.

\begin{prop} \label{prop:unpuncturedmesh}
Let $i,j,k,l$ be four boundary vertices of $S$ in clockwise order around the
boundary of $S$. Then:
$$n_{ik}n_{jl}=n_{ij}n_{kl}+n_{li}n_{jk}.$$
\end{prop}

Finally, we note that the numbers $(i,j)$ can be interpreted as
specialisations of cluster variables by work of Fomin and
Zelevinsky~\cite[12.2]{fominzelevinsky03}. Let $\mathcal{A}$ be a cluster
algebra of type $A_{N-3}$ with trivial coefficients. Then the cluster
variables of $\mathcal{A}$ are in bijection with the set consisting of
all of the diagonals of $S$ (using also~\cite{fominzelevinsky01}).
For distinct vertices $i,j$, let $x_{ij}$ be the cluster variable
corresponding to the diagonal $D_{ij}$ of $S$.
By~\cite[3.1]{fominzelevinsky02}, $x_{ij}$ is
a Laurent polynomial in the cluster variables corresponding to the diagonals
of $\T$.

The following result is proved in~\cite[\S 3]{propp} (for the $n_{ij}$),
but we include a proof for the convenience of the reader.
We also note that a connection between frieze patterns and cluster
algebras was first explicitly given in~\cite{calderochapoton06}. 

\begin{theorem} \label{thm:clustertypeA}
Let $i,j$ be distinct vertices on the boundary of $S$. Then $(i,j)$ is
equal to the number $u_{ij}$ obtained from $x_{ij}$ by specialising each
of the cluster variables corresponding to the arcs of $\T$ to $1$.
\end{theorem}

\begin{proof}
This follows from~\cite{fominzelevinsky03} together with
Theorem~\ref{thm:bracketequalsn} and Proposition~\ref{prop:unpuncturedmesh},
which show that the $(i,j)$ and the $u_{ij}$ both satisfy the
relations of Proposition~\ref{prop:unpuncturedmesh}.
Since $(i,j)=1$ is equal to the specialisation of $x_{ij}$ whenever $D_{ij}$
is an arc in $\T$, the equality for arbitrary diagonals follows from
iterated application of these relations.
\end{proof}

%
\section{Construction of frieze patterns of type $D_N$}
\label{sec:mesh-rel}
%

Our aim in this section is to prove the main result, that the numbers
$m_{ii}$, $m_{ij}$ and $m_{i0}$ form a frieze
pattern of type $D$. In order to do that we have to prove that the
frieze relations (see Section~\ref{sec:introduction}) hold. We now work
with the disc $S_{\odot}$ with puncture $0$ and $N$ marked points on the
boundary labelled $\{1,2,\ldots ,N\}$.
We first need the following result.

\begin{lemma} \label{lem:onecase}
Let $\T$ be a triangulation of $S_{\odot}$, and let
$i, j\in\{1,\dots,N\}$ with $j\neq i+1$.

(1) If $D_{ij}\in\T$ then $m_{ij}=1$.

(2) If $D_{i0}\in\T$ then $m_{i0}=1$.
\end{lemma}

\begin{proof}
(1) Consider the truncated disc $S_{ij}$ with boundary
$B_{ij}$ and $D_{ij}$. Since $D_{ij}\in\T$, none of the triangles
of $\T$ are split by $D_{ij}$ and $\T|S_{ij}$ is a triangulation
of $S_{ij}$. Suppose first that $i \neq j$.
Then $m_{ij}$ is the number of matchings between
$\{i+1,\dots,j-1\}$ and $\T|_{S_{ij}}$. Such a matching is
a matching between $\T|_{S_{ij}}$ and all of the vertices of
the unpunctured disc $S_{ij}$ except $i$ and $j$ which
are neighbours on the boundary of $S_{ij}$. By
Remark~\ref{rem:unpuncturedone}, $m_{ij}=1$.

Suppose instead that $i=j$.
Then $m_{ii}$ is the number of matchings between
$\{1,2,\ldots ,N\}\setminus \{i\}\cup \{0\}$.
Such a matching is a matching between $\T|_{S_{ii}}$ and all of
the vertices of the unpunctured disc $P(i)$ except
$i$ and $i'$ which are joined by the arc $D_{ii'}$ in $P(i)$
induced by the arc $D_{ii}$ of $\T$. By Remark~\ref{rem:unpuncturedone},
$m_{ii}=1$.

(2) If $D_{i0}\in\T$ then $i$ and $0$ are neighbours in
$P(i)$ and the statement follows with the same reasoning. 
\end{proof}

\begin{remark} \label{rem:counting}
We recall that in a triangulation of an unpunctured disc with $N$
marked points on the boundary, there are $N-2$ triangles. It follows that
there are no matchings between a triangulation of the disc and
a subset of the boundary vertices of cardinality greater than $N-2$.
\end{remark}

\begin{definition} 
Let $M$ be a matching between a subset $I$ of $\{1,\ldots ,N\}\cup \{0\}$
and a triangulation $\T$ of $S_{\odot}$. Let $P$ be a subset
of $S_{\odot}$ consisting of a union of triangles of $\T$.
We denote by $M|_P$ the restriction of $M$ to $P$: this is the
matching between $\T|_P$ and $I\cap P$ obtained by allocating each
vertex its corresponding triangle in $M$ whenever that triangle is
contained in $P$.
\end{definition}

In order to do some computations of numbers of matchings, we need the
following simple observation:

\begin{lemma} \label{lem:decomposition}
Let $\T$ be a triangulation of $S_{\odot}$, and let
$S_{\odot}=P\cup Q$ be a decomposition of $S_{\odot}$ into two subsets
with common boundary given by arcs of $\T$. Given a subset
$I$ of $\{1,2,\ldots ,N\}\cup \{0\}$, let $J$ denote the subset of $I$
consisting of vertices on the common boundary between $P$ and $Q$.
Then the number of matchings between $\T$ and $I$ is given by
$$|\M (I,\T)|=\sum_{J',J''\,:J=J'\sqcup J''} n_{J',J''}.$$
Here $n_{J',J''}$ is the number of matchings between $\T$ and $I$
in which all vertices of $J'$ are allocated triangles in $P$
and all vertices of $J''$ are allocated triangles in $Q$. It is
given by the product of the number of matchings between
$\T|_P$ and $I\setminus J'$ and the number of matchings between
$\T|_Q$ and $I\setminus J''$.
\end{lemma}

\begin{proof}
Let $M$ be a matching between $I$ and $\T$.
Given a pair $M_P,M_Q$ of matchings for $\T|_P$ and $\T|_Q$ which are
compatible in the sense that any vertex of $J$ is allocated precisely one
triangle in either $\T|_P$ or $\T|_Q$ (but not both), we
can put them together to obtain a matching $M$ between $I$ and $\T$.
This gives us a bijection between matchings $M$ between $I$ and $\T$
and compatible pairs $M',M''$ of matchings for $\T|_P$ and
$\T|_Q$. The result follows from dividing up such pairs according to
the allocation of the vertices of $J$ to triangles in $P$ or in $Q$.
\end{proof}

\begin{lemma} \label{lem:degeneratecase}
Let $\T$ be a triangulation of $S_{\odot}$.
Let $d_0$ denote the number of triangles of $\T$ incident with the puncture,
$0$, and let $i$ denote any boundary vertex of $S_{\odot}$.
Then we have $d_0m_{i0}=m_{ii}$.
\end{lemma}
\begin{proof}
If there are at least two arcs incident with $0$ in $\T$, let $j$ be the
first boundary vertex strictly clockwise from $i$ such that $D_{j0}\in \T$,
and let $k$ be the first boundary vertex anticlockwise from $i$
(possibly equal to $i$) such that $D_{k0}\in \T$.
Since there are no arcs of the form $D_{l0}\in \T$ for $l\in V_{k+1,j-1}$,
we have $D_{kj}\in \T$.
If there is only one arc $D_{l0}$ incident with $0$ in $\T$ then we
take $j=k=l$. In this case $D_{jj}\in \T$.

If there is more than one arc incident with $0$ in $\T$, then we have a
decomposition of the disc $S_{\odot}$ into smaller discs $P(i)$ and $Q(i)$
(Definition~\ref{def:piqi}).
We shall apply Lemma~\ref{lem:decomposition} in these cases.
For $i,j$ any pair of boundary vertices of $P(i)$, let $p_{ij}$ denote
the number of matchings of $\T|_{P(i)}$ with all of the boundary vertices of
$P(i)$ except $i$ and $j$. Let $q_{ij}$ denote the corresponding number for
a pair of boundary vertices of $Q(i)$.

\noindent {\bf Case (I):}
We assume first that $i,j$ and $k$ are all distinct.
See Figure~\ref{fig:degeneratecaseI}.
If, in a matching in $\M_{ii}$, both $j$ and $k$ are allocated triangles in
$P(i)$ then $0$ must be allocated a triangle in $Q(i)$ by
Remark~\ref{rem:counting} applied to $P(i)$.
By restriction of such a matching to $P(i)$ we obtain a matching
for $\T|_{P(i)}$ in which $i$ and $0$ are not allocated triangles
and by restriction to $Q(i)$ we obtain a matching for $\T|_{Q(i)}$ in
which $j$ and $k$ are not allocated triangles.
Thus there are $p_{i0}q_{jk}$ matchings in $\M_{ii}$ in which
$j$ and $k$ are both allocated triangles in $P(i)$.
Using the fact that $p_{i0}=m_{i0}$
(by definition) and Theorem~\ref{thm:broline}, we see that this is equal to
$m_{i0}(d_0-1)$.
If $j$ is allocated a triangle
in $Q(i)$ and $k$ is allocated a triangle in $P(i)$ then $0$ is allocated
a triangle in $P(i)$ by Remark~\ref{rem:counting} applied to $Q(i)$, and we
obtain
$p_{ij}q_{0k}=p_{ij}\cdot 1=p_{ij}$ matchings of this type by
Remark~\ref{rem:unpuncturedone}.
Similarly, there are $p_{ij}$ matchings in $\M_{ii}$ in which
$j$ is allocated a triangle in $P(i)$ and $k$ is allocated a triangle in
$Q(i)$.

There are no matchings in $\M_{ii}$ in which both $j$ and $k$ are both
allocated triangles in $Q(i)$, by Remark~\ref{rem:counting} applied to
$Q(i)$, so we have covered all cases.
By Lemma~\ref{lem:decomposition} we have that
$m_{ii}=|\M_{ii}|=(d_0-1)m_{i0}+p_{ij}+p_{ik}$.
By Proposition~\ref{prop:unpuncturedmesh},
$p_{i0}p_{jk}=p_{ij}p_{0k}+p_{ik}p_{0j}$. Since $D_{kj}\in \T$,
$p_{jk}=1$ by Remark~\ref{rem:unpuncturedone}, and we also have
by Remark~\ref{rem:unpuncturedone} that $p_{0j}=p_{0k}=1$, so
$p_{ij}+p_{ik}=p_{i0}=m_{i0}$. We obtain
$$m_{ii}=(d_0-1)m_{i0}+m_{i0}=d_0m_{i0}$$
as required.

\begin{figure}\label{fig:degeneratecaseIandII}
\begin{center}
\subfigure[Case I]
   {\includegraphics[scale=.4]{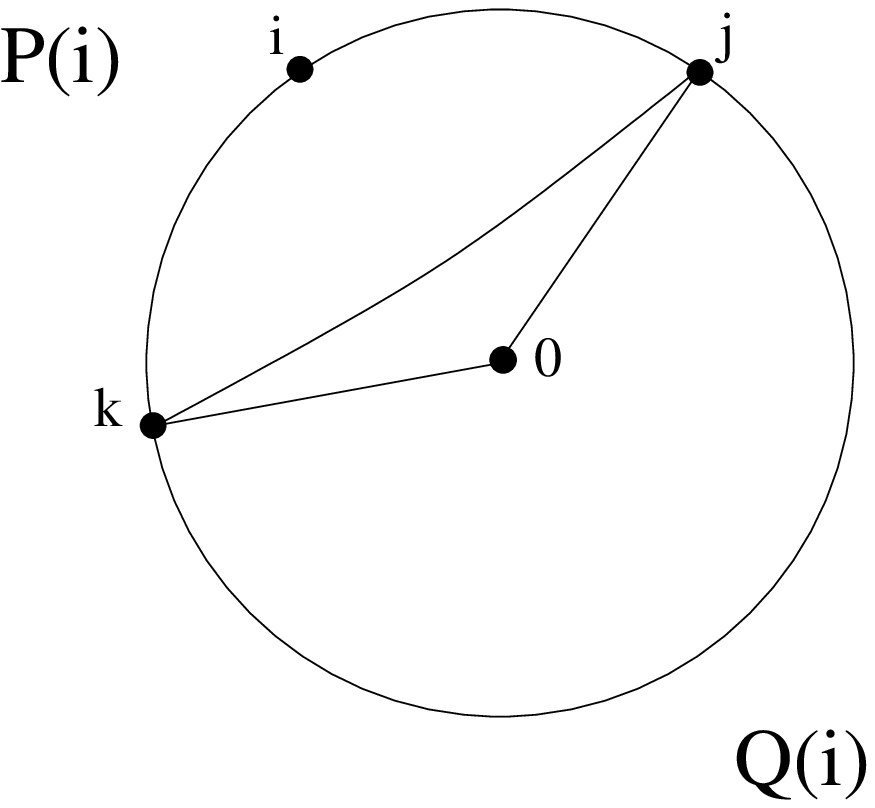}\label{fig:degeneratecaseI}}
\hskip 2cm
\subfigure[Case II]
   {\includegraphics[scale=.4]{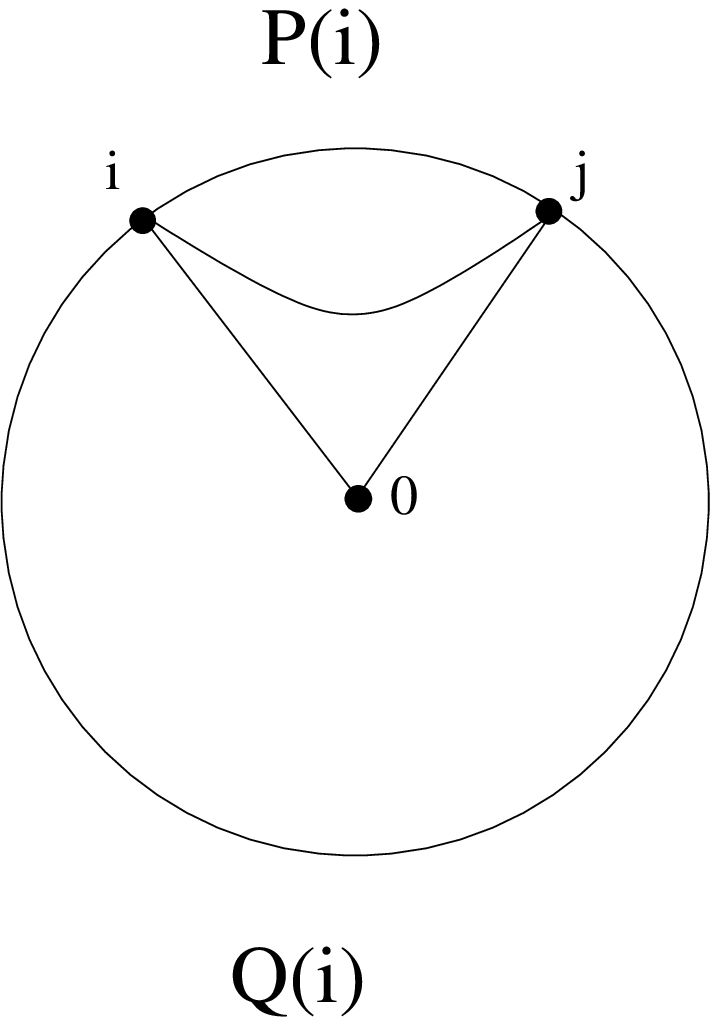}\label{fig:degeneratecaseII}}
\end{center}
\caption{Proof of Lemma~\ref{lem:degeneratecase}, Cases (I) and (II)}
\end{figure}

\noindent {\bf Case (II):}
We assume next that $i=k$.
See Figure~\ref{fig:degeneratecaseII}.
If, in a matching in $\M_{ii}$, $j$ is allocated a triangle in $P(i)$,
then $0$ must be allocated a triangle in $Q(i)$ by
Remark~\ref{rem:counting} applied to $P(i)$, and we see that
there are $p_{i0}q_{ij}$ matchings of this type. By
Theorem~\ref{thm:broline} this is equal to $m_{i0}(d_0-1)=d_0-1$,
noting that $m_{i0}=1$ by Lemma~\ref{lem:onecase}.
If $j$ is allocated a triangle in $Q(i)$, then $0$ must be allocated
a triangle in $P(i)$ by Remark~\ref{rem:counting} applied to $Q(i)$.
We obtain $p_{ij}q_{0i}=1\cdot 1=1$ matchings of this type by
Remark~\ref{rem:unpuncturedone}.
By Lemma~\ref{lem:decomposition} we obtain a total number of matchings
$m_{ii}=|\M_{ii}|=d_0-1+1=d_0=m_{i0}d_0$ as required.

\noindent {\bf Case (III):}
We assume next that $i\not=j=k$.
In a matching in
$\M_{ii}$, only one of $j,j'$ is allocated a triangle.
See Figure~\ref{fig:degeneratecaseIII}.

We see that $m_{ii}$ is given by the the sum of the number of matchings
between $\T|_{P(i)}$ and $V_{i+1,i-1}\cup \{0\}$ and the number of matchings
between $\T|_{P(i)}$ and $V_{i+1,i-1}\cup \{0\}$ with $j$ replaced by $j'$.
There are $p_{ij}$ matchings of the first kind and
$p_{ij'}$ matchings of the second kind, making a total of $p_{ij}+p_{ij'}$.
We have that $p_{i0}p_{jj'}=p_{ij}p_{0j}+p_{ij'}p_{0j'}$,
by Proposition~\ref{prop:unpuncturedmesh}, so $m_{i0}=p_{i0}=p_{ij}+p_{ij'}$
by Remark~\ref{rem:unpuncturedone},
noting that $0,j$ and $0,j'$ are adjacent on the boundary of $P(i)$
and $D_{j'j}$ lies in $\T|_{P(i)}$.
We obtain $m_{ii}=m_{i0}=m_{i0}d_0$ as required.

\begin{figure}
\begin{center}
\includegraphics[scale=.4]{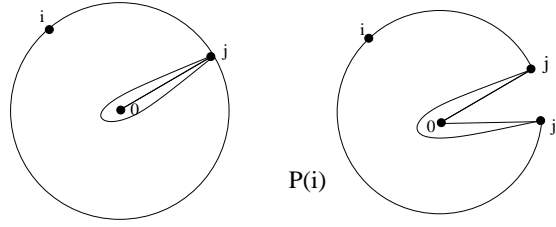}
\caption{Proof of Lemma~\ref{lem:degeneratecase}, Case (III).}
\label{fig:degeneratecaseIII}
\end{center}
\end{figure}

\noindent {\bf Case (IV)}
We finally assume that $i=j=k$, so that there is a unique arc incident
with $0$ in $\T$ given by $D_{i0}$.
See Figure~\ref{fig:degeneratecaseIV}.
We have that $m_{ii}$ is given by the
the number of matchings between $V_{i+1,i-1}\cup \{0\}$ and $\T|_{P(i)}$,
so we obtain $m_{ii}=p_{ii'}=1$ by Remark~\ref{rem:unpuncturedone}
since $D_{i'i}\in \T|_{P(i)}$.
Hence $m_{ii}=d_0m_{i0}$ as required, since $d_0=1$ and $m_{i0}=1$
by Lemma~\ref{lem:onecase}.

\begin{figure}
\begin{center}
\includegraphics[scale=.4]{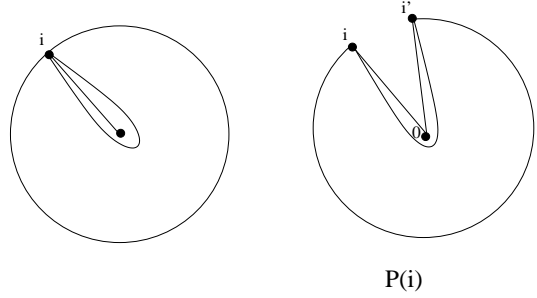}
\caption{Proof of Lemma~\ref{lem:degeneratecase}, Case (IV).}
\label{fig:degeneratecaseIV}
\end{center}
\end{figure}

The proof of Lemma~\ref{lem:degeneratecase} is complete.
\end{proof}

In order to prove the frieze relations hold, we first need the following:

\begin{lemma} \label{lem:extension}
Let $\T$ be a triangulation of $S_{\odot}$.
Let $i,j$ be boundary vertices of $S_{\odot}$ with $j\neq i+1$.
Then the restriction $\T|_{S_{ij}}$ can be extended to a
triangulation of the entire disc $S$ with different marked points and
the puncture removed.
\end{lemma}

\begin{figure}
\begin{center}
\includegraphics[scale=.7]{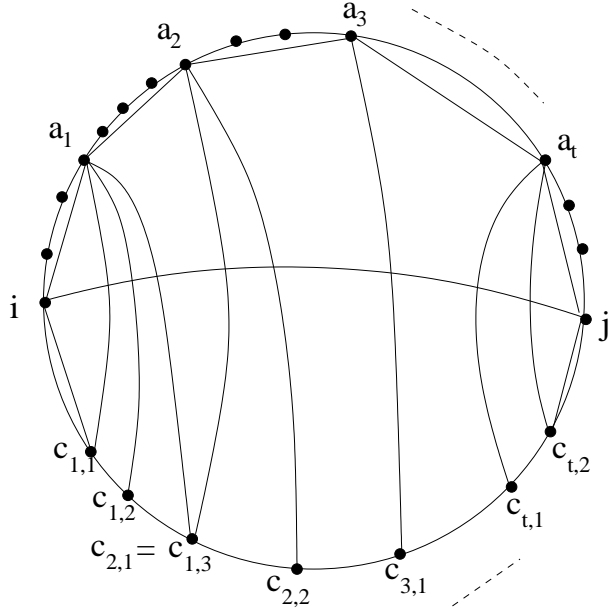}
\caption{Completing to a triangulation of the unpunctured disc.}
\label{fig:completion}
\end{center}
\end{figure}

\begin{proof}
Let $a_1,a_2,\ldots ,a_t$ be the vertices on the arc $B_{i,j}$
(distinct from $i$ and $j$) in order clockwise from $i$, which lie on
arcs of $\T$ whose other end lies in $V_{j+1,i-1}$
For each $l$, let $b_{l,1},b_{l,2},\ldots ,b_{l,u_l}$ be the end-points
(other than $a_l$) of the arcs of $\T$ incident with $a_l$ whose
other end lies in $V_{j+1,i-1}$.

The arc $D_{i,a_1}$ must lie in $\T$, else the region on the $i$-side
of the arc $D_{a_1,b_{l,1}}$ will not be a triangle, since the other
end-point of any arc incident with any of the vertices
on the boundary arc $B_{i+1,a_1-1}$ can only be a vertex on the boundary
arc $B_{i,a_1}$. Similarly, the arcs $D_{a_m,a_{m+1}}$ must all be in $\T$,
as must $D_{a_t,j}$.

Remove the vertices $V_{j+1,i-1}$ and introduce new vertices on the boundary
between $j$ and $i$ labelled
$$c_{1,1},c_{1,2},\ldots ,c_{1,u_1},c_{2,1},\ldots  ,c_{l,u_l},$$
going anticlockwise from $i$ to $j$. We identify $c_{m,u_m}$ with
$c_{m+1,1}$ for each $m$. We replace the part of each arc
$D_{a_m,b_{l,m}}$ below the arc $D_{i,j}$
with a new arc linking $c_{l,m}$ with the intersection of $D_{a_m,b_{l,m}}$ and
$D_{i,j}$, and add arcs $D_{c_{1,1},i}$ and $D_{j,c_{t,u_t}}$.
See Figure~\ref{fig:completion} for an example.
In this way we can complete $\T|_{S_{i,j}}$ to a triangulation
of the disc $S$, with new marked points on the boundary and the
puncture removed, as required.
\end{proof}

We first note the following consequence:

\begin{lemma} \label{lem:positivenumbers}
Let $\T$ be a triangulation of $S_{\odot}$.
Let $i,j$ be boundary vertices of $S_{\odot}$ with $j\neq i,i+1$.
Then the numbers $m_{ii}$, $m_{ij}$ and $m_{i0}$
as in Definition~\ref{def:matchingnumbers} are all positive.
\end{lemma}

\begin{proof}
By Definition~\ref{def:matchingnumbers}, the $m_{i0}$ can be interpreted
as entries in Conway-Coxeter frieze patterns (using
Theorem~\ref{thm:bracketequalsn}), so are positive. By Lemma~\ref{lem:extension}
the same argument applies to the $m_{ij}$. Then the $m_{ii}$ are positive
by Lemma~\ref{lem:degeneratecase}.
\end{proof}

We can now prove that the frieze relations hold:

\begin{prop} \label{prop:relations}
Let $\T$ be a triangulation of $S_{\odot}$.
Let $i,j$ be boundary vertices of $S_{\odot}$ with $j\neq i,i+1$,
Let the numbers $m_{ii}$, $m_{ij}$ and $m_{i0}$ be as in
Definition~\ref{def:matchingnumbers}. Then the following hold:
\begin{eqnarray}
m_{ij}\cdot m_{i+1,j+1} & = & m_{i+1,j}\cdot m_{i,j+1} + 1,
\mbox{ provided $j\not=i$ or $i+1$}
\label{eqn:relation1prop}
\\
m_{i,i-1}\cdot m_{i+1,i} & = & m_{i+1,i-1}\cdot m_{ii}\cdot m_{i0}+1
\label{eqn:relation2prop} \\
m_{ii}\cdot m_{i+1,0} & = & m_{i+1,i} + 1
\label{eqn:relation3prop} \\
m_{i0}\cdot m_{i+1,i+1} & = & m_{i+1,i} + 1
\label{eqn:relation4prop}
\end{eqnarray}
\end{prop}

\begin{proof}
We prove each relation in turn. \\
\noindent {\bf Proof of (1):}
By Lemma~\ref{lem:extension} we can complete $\T|_{S_{i,j+1}}$ to a triangulation
$\T'$ of the disc $S$ (with the puncture removed) and new marked points
on the boundary.
We see that $m_{i,j+1}$ is the number of matchings between $\T'$
and $V_{i+1,j-1}$. It is clear that
$m_{ij}$, $m_{i+1,j}$ and $m_{i,j+1}$ are the numbers of matchings between
$\T'$ and $V_{i+1,j-1}$, $V_{i+2,j-1}$ and $V_{i+1,j}$
respectively. Hence the frieze relation (1) holds by
Proposition~\ref{prop:unpuncturedmesh}.

\noindent {\bf Proof of (2):}
We adopt the same notation for $j$ and $k$ as in the proof of
Lemma~\ref{lem:degeneratecase}.

Let $\widetilde{\M}_{ii}$ denote the set of matchings between
$I:=\{1,2,\ldots ,N\}\setminus \{i\}$ and $\T|_{S_{ii}}$.
Set $\widetilde{m}_{ii}=|\widetilde{\M}_{ii}|$.
It is straightforward to show (along the lines of the proof of (1))
that, for any boundary vertex $i$ of $S_{\odot}$,
$$m_{i,i-1}m_{i+1,i}=m_{i+1,i-1}\widetilde{m}_{ii}+1.$$
Hence for (2) it is enough to show that $\widetilde{m}_{ii}=m_{i0}m_{ii}$.
By Lemma~\ref{lem:degeneratecase}, this is equivalent to showing that
$\widetilde{m}_{ii}=d_0m_{i0}^2$.
 
If there is more than one arc incident with $0$ in $\T$, then we have a
decomposition of the disc $S_{\odot}$ into smaller discs $P(i)$ and $Q(i)$
(see Definition~\ref{def:piqi}).
We shall apply Lemma~\ref{lem:decomposition} in these cases. Let $p_{ij}$
and $q_{ij}$ be as in Lemma~\ref{lem:degeneratecase}. 

\noindent {\bf Proof of (2), Case (I):}
We assume first that $i,j$ and $k$ are distinct.
See Figure~\ref{fig:part2caseI}.

\begin{figure}\label{fig:part2caseIandII}
\begin{center}
\subfigure[Case I]
   {\includegraphics[scale=.4]{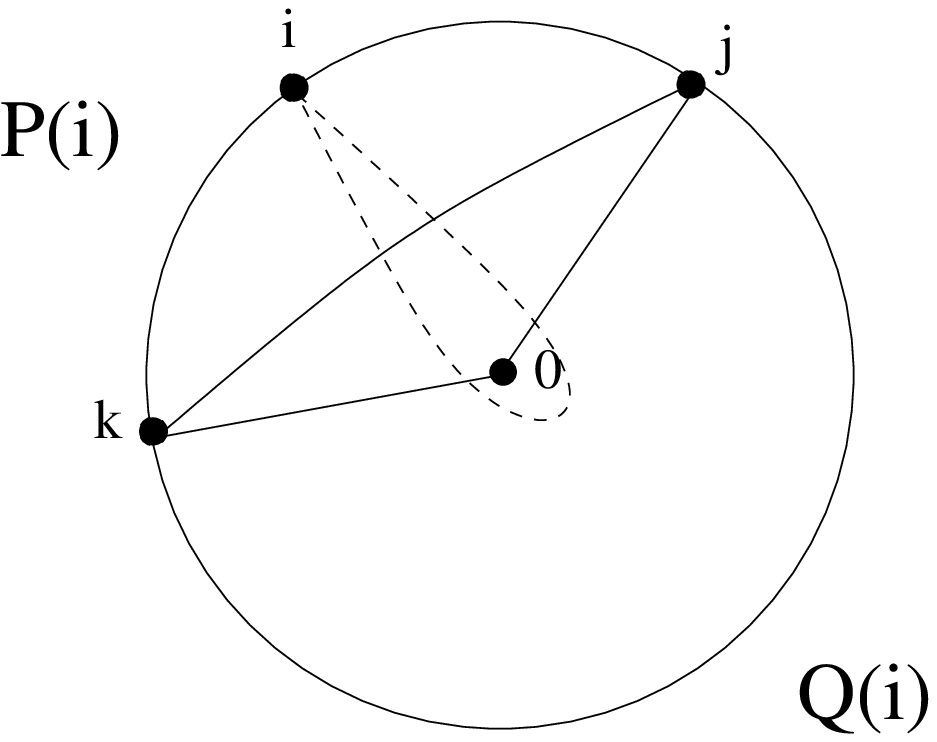}\label{fig:part2caseI}}
\hskip 2cm
\subfigure[Case II]
   {\includegraphics[scale=.4]{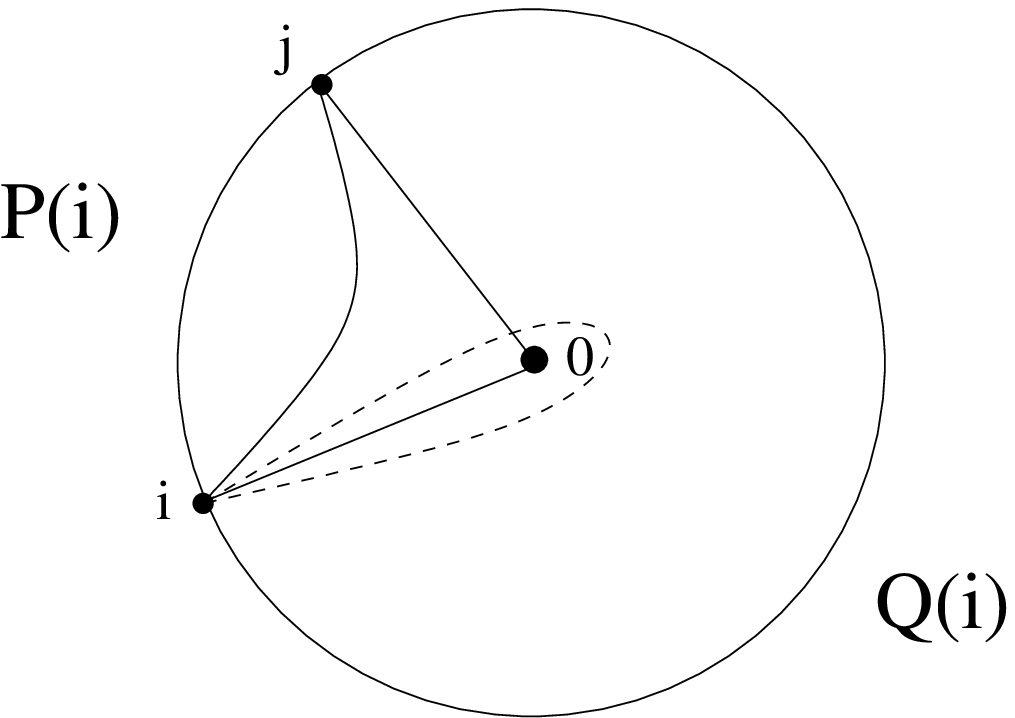}\label{fig:part2caseII}}
\end{center}
\caption{Proof of Proposition~\ref{prop:relations}(2), Cases (I) and (II)}
\end{figure}

(a) Suppose first that in a matching in $\widetilde{\M}_{ii}$,
$j,k$ are both allocated triangles in $P(i)$. Then
restricting the matching to $P(i)$ we obtain a matching between
$\T|_{S_{ii}}$ and $V_{i+1,j}\cup V_{k,i-1}$.
Since the arc $D_{ii}$ splits $P(i)$ completely, there are
$|\M(V_{i+1,j},\T|_{P(i)})|\cdot |\M(V_{k,i-1},\T|_{P(i)})|$ of these.
By Theorem~\ref{thm:broline}, this is equal to $m_{i0}^2$. Restricting $\T$ to
$Q(i)$, no triangles are split by $D_{ii}$. Therefore there are
$|\M(V_{j+1,k-1},\T|_{Q(i)})|$ matchings of this type in $\widetilde{\M}_ii$.
By Theorem~\ref{thm:broline}, this is equal to $\M(\{0\},\T|_{Q(i)})=d_0-1$.
We see that there are $(d_0-1)m_{i0}^2$
matchings in $\widetilde{\M}_{ii}$ in which $j$ and $k$ are both allocated
triangles in $P(i)$.

(b) Arguing as above we see that there are
$|\M(V_{k+1,i-1},\T|_{P(i)})|\cdot |\M(V_{i+1,j},\T|_{P(i)})|$
possible restrictions to $P(i)$ of a matching in
$\widetilde{\M}_{ii}$ in which $j$ is allocated a triangle in $P(i)$ and $k$ a
triangle in $Q(i)$. By Theorem~\ref{thm:broline} this is equal to
$p_{ki}p_{i0}=p_{ki}m_{i0}$.
Similarly we see that there are $|\M(V_{j+1,k},\T|_{Q(i)})|=q_{j0}=1$
possible restrictions to $Q(i)$, the last equality using
Remark~\ref{rem:unpuncturedone}. We get a total of $p_{ki}m_{i0}$
matchings of this type in $\widetilde{\M}_{ii}$.

(c) Arguing as in (b) we see that there are $p_{ji}m_{i0}$ matchings in
$\widetilde{\M}_{ii}$ in which $j$ is allocated a triangle in $Q(i)$ and
$k$ a triangle in $P(i)$.

(d) By Remark~\ref{rem:counting} for $Q(i)$, $j$ and $k$ cannot both be
allocated triangles in $Q(i)$ for any matching in $\widetilde{\M}_{ii}$.

By Lemma~\ref{lem:decomposition} we see that
$\widetilde{m}_{ii}=(d_0-1)m_{i0}^2+(p_{ki}+p_{ij})p_{i0}$. But
by Proposition~\ref{prop:unpuncturedmesh}, $p_{i0}=p_{ij}+p_{ki}$, so
$\widetilde{m}_{ii}=d_0m_{i0}^2$ as required.

\noindent {\bf Proof of (2), Case (II):}
We next assume that $i=k\not=j$.
See Figure~\ref{fig:part2caseII}.

(a) There are $|\M(V_{i+1,j},\T|_{P(i)})|$ possible restrictions to $P(i)$
of a matching in $\widetilde{\M}_{ii}$ in which $j$ is allocated a vertex
in $P(i)$. There are $|\M(V_{j+1,i-1},\T|_{Q(i)})|$ possible
restrictions to $Q(i)$, giving a total of
$$|\M(V_{i+1,j},\T|_{P(i)})|\cdot |\M(V_{j+1,i-1},\T|_{Q(i)})|=p_{i0}
|\M(\{0\},\T|_{Q(i)})|=(d_0-1)$$
by Theorem~\ref{thm:broline} and the fact that $p_{i0}=m_{i0}=1$
(by Lemma~\ref{lem:onecase}).

(b)There are
$|\M(V_{i+1,j-1},\T|_{P(i)})|$ possible restrictions to $P(i)$ of a matching
in $\widetilde{\M}_{ii}$ in which $j$ is allocated a triangle in $Q(i)$.
There are $|\M(V_{j,i-1},\T|_{Q(i)})|$ possible restrictions to $Q(i)$.
By Theorem~\ref{thm:broline} this is equal to $p_{ij}q_{i0}$,
which is equal to $1$ by Lemma~\ref{lem:onecase}.

By Lemma~\ref{lem:decomposition}, we have $\widetilde{m}_{ii}=d_0-1+1=d_0$.
This equals $m_{i0}d_0$ by Lemma~\ref{lem:onecase} so we are done.

\noindent {\bf Proof of (2), Case (III):}
We next assume that $j=k\not=i$.
See Figure~\ref{fig:part2caseIII}.
\begin{figure}
\begin{center}
\includegraphics[scale=.4]{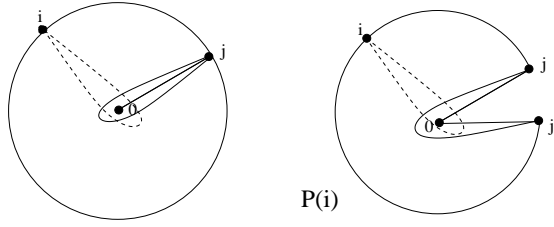}
\caption{Proof of Proposition~\ref{prop:relations}(2), Case (III).}
\label{fig:part2caseIII}
\end{center}
\end{figure}

A matching in $\widetilde{\M}_{ii}$ induces a matching of $\T|_{P(i)}$
in which either $j$ or $j'$ is allocated a triangle but not both.
Since $P(i)$ is split by $D_{ii}$ completely, there are
$|\M(V_{j',i-1},\T|_{P(i)})|\cdot |\M(V_{i+1,j-1},\T|_{P(i)})|$
matchings of the first kind. By Theorem~\ref{thm:broline} this is equal to
$p_{i0}p_{ij}$.
Similarly there are $p_{ij'}p_{i0}$ matchings of the second kind.

We get a total of $p_{i0}(p_{ij}+p_{ij'})$ matchings in $\widetilde{\M}_{ii}$.
By Proposition~\ref{prop:unpuncturedmesh} we see that $p_{i0}=p_{ij}+p_{ij'}$,
and we obtain $\widetilde{m}_{ii}=p_{i0}^2=m_{i0}^2=d_0m_{i0}^2$ as required,
since $d_0=1$.

\noindent {\bf Proof of (2), Case (IV):}
Finally we consider the case where $i=j=k$.
See Figure~\ref{fig:part2caseIV}.
\begin{figure}
\begin{center}
\includegraphics[scale=.4]{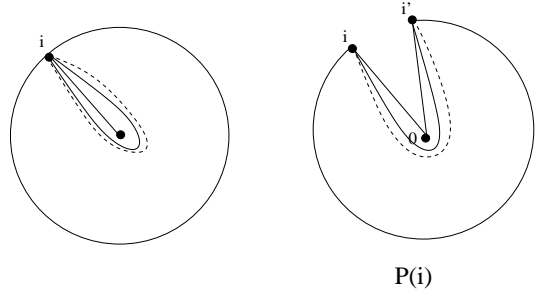}
\caption{Proof of Proposition~\ref{prop:relations}(2), Case (IV).}
\label{fig:part2caseIV}
\end{center}
\end{figure}

A matching in $\widetilde{\M}_{ii}$ induces a matching of $\T|_{P(i)}$
in which all boundary vertices except $i$ and $i'$ are allocated a
triangle. Since $D_{i'i}$ is an arc in $\T|_{P(i)}$, we see by
Remark~\ref{rem:unpuncturedone} that there is only one possible such matching.
So $\widetilde{m}_{ii}=1=d_0m_{i0}$, since $m_{i0}=1$
by Lemma~\ref{lem:onecase}.

\noindent {\bf Proof of (3):}
We note that by Lemma~\ref{lem:degeneratecase}, it is sufficient to show
that $d_0m_{i0}m_{i+1,0}=m_{i+1,i}+1$.

\noindent {\bf Proof of (3), Case (I):}
We first assume that $i,i+1,j$ and $k$ are all distinct.
See Figure~\ref{fig:part3caseI}.

\begin{figure}\label{fig:part3caseIandII}
\begin{center}
\subfigure[Case I]
   {\includegraphics[scale=.4]{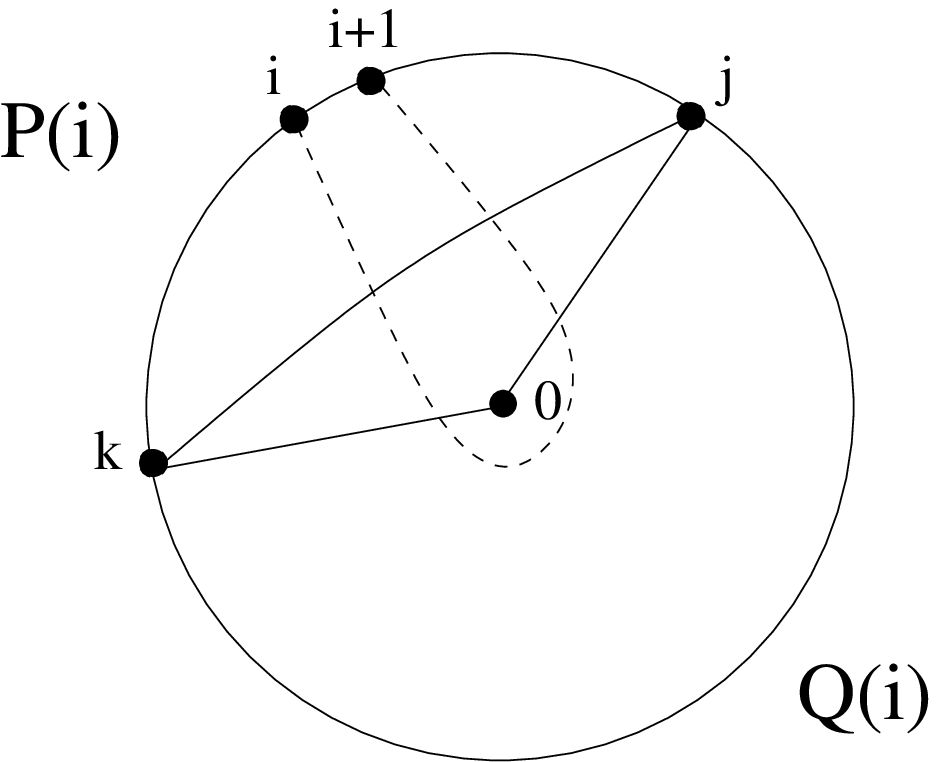}\label{fig:part3caseI}}
\hskip 2cm
\subfigure[Case II]
   {\includegraphics[scale=.4]{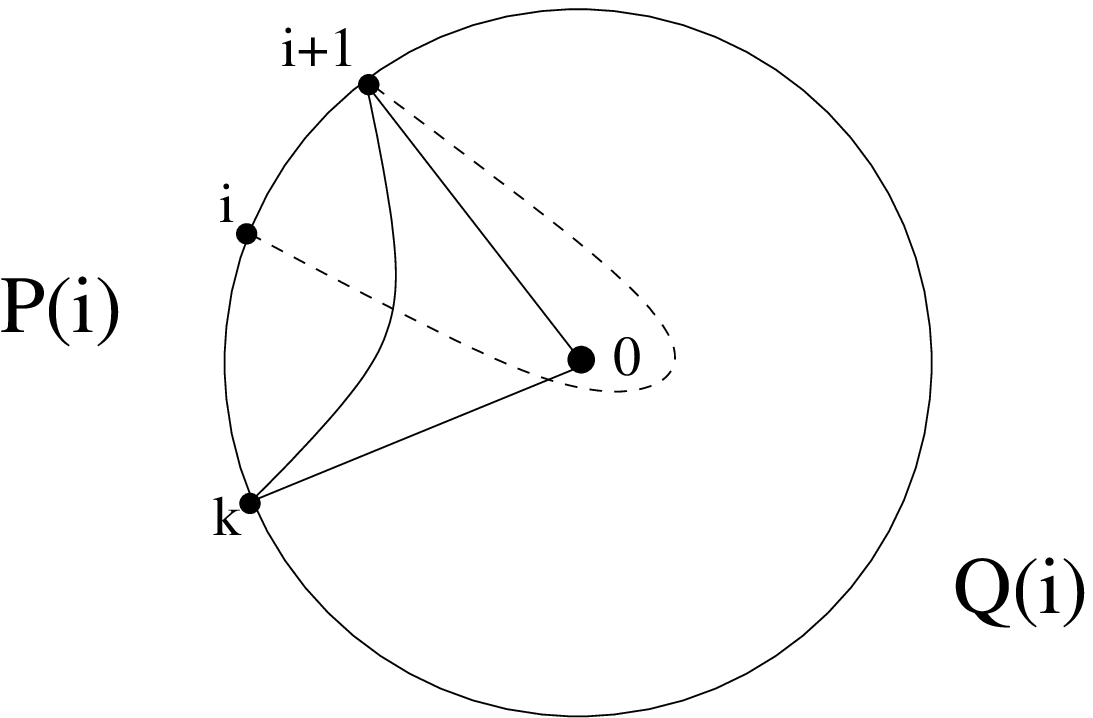}\label{fig:part3caseII}}
\end{center}
\caption{Proof of Proposition~\ref{prop:relations}(3), Cases (I) and (II)}
\end{figure}

(a) There are $|\M(V_{i+2,j},\T|_{P(i)})|\cdot |\M(V_{k,i-1},\T|_{P(i)})|$
possible restrictions to $P(i)$ of a matching in $\M_{i+1,i}$
in which both $j$ and $k$ are allocated triangles in $P(i)$,
since $D_{i+1,i}$ splits $P(i)$ completely.
By Theorem~\ref{thm:broline} (and the definition of $m_{i0}$ and $m_{i+1,0}$),
this equals $m_{i+1,0}m_{i0}$.
There are $|M(V_{j+1,k-1},\T|_{Q(i)})|$
possible restrictions to $Q(i)$, which by
Theorem~\ref{thm:broline} is equal to
$|M(\{0\},\T|_{Q(i)})|=d_0-1$.
Hence there are $(d_0-1)m_{i0}m_{i+1,0}$ matchings in $\M_{i+1,i}$
in total in which both $j$ and $k$ are allocated triangles in $P(i)$.

(b) There are
$|\M(V_{i+2,j},\T|_{P(i)})|\cdot |\M(V_{k+1,i-1},\T|_{P(i)})|$
possible restrictions to $P(i)$ of a matching in $\M_{i+1,i}$ in
which $j$ is allocated a triangle in $P(i)$ and $k$ is allocated
a triangle in $Q(i)$. By Theorem~$\ref{thm:broline}$ this is equal to
$p_{i+1,0}p_{ik}=m_{i+1,0}p_{ik}$. There are
$|M(V_{j+1,k},\T|_{Q(i)})|$
possible restrictions to $Q(i)$, which is equal to $q_{0j}$ by
Theorem~\ref{thm:broline} and thus equal to $1$ by
Remark~\ref{rem:unpuncturedone}.
Thus we see that there are a total
of $p_{ik}m_{i+1,0}$ matchings in $\M_{i+1,i}$ in which $j$ is
allocated a triangle in $P(i)$ and $k$ is allocated a triangle in $Q(i)$.

(c) Arguing in a similar way to (b) we see that there are $p_{i+1,j}m_{i0}$
matchings in $\M_{i+1,i}$ in which $j$ is allocated a triangle in
$Q(i)$ and $k$ is allocated a triangle in $P(i)$. 

(d) We note that it is not possible for both $j$ and $k$ to be allocated
triangles in $Q(i)$ in a matching in $\M_{i+1,i}$,
by Remark~\ref{rem:counting} applied to $Q(i)$.

By Lemma~\ref{lem:decomposition} we see that
$$m_{i+1,i}=(d_0-1)m_{i0}m_{i+1,0}+p_{ik}m_{i+1,0}+p_{i+1,j}m_{i0}.$$
By Proposition~\ref{prop:unpuncturedmesh} we have that
$p_{i0}p_{jk}=p_{ij}p_{0k}+
p_{ik}p_{0j}$, so $p_{i0}=p_{ij}+p_{ik}$, noting that $D_{kj}$ is an arc
in $\T|_{P(i)}$, so $p_{jk}=1$. Hence
\begin{eqnarray*}
m_{i+1,i} &=& 
(d_0-1)m_{i0}m_{i+1,0}+ m_{i+1,0}(p_{i0}-p_{ij})+m_{i0}p_{i+1,j} \\
&=&
(d_0-1)m_{i0}m_{i+1,0}+ m_{i+1,0}(m_{i0}-p_{ij})+m_{i0}p_{i+1,j} \\
&=&
(d_0-1)m_{i0}m_{i+1,0}+ m_{i0}m_{i+1,0}-m_{i+1,0}p_{ij}+m_{i0}p_{i+1,j} \\
&=&
(d_0-1)m_{i0}m_{i+1,0}+ m_{i0}m_{i+1,0}-1 \\
&=&
d_0m_{i0}m_{i+1,0}-1,
\end{eqnarray*}
using the fact that
\begin{eqnarray*}
p_{i+1,0}p_{ij} &=& p_{i+1,j}p_{i0}+p_{0j}p_{i,i+1} \\
&=& p_{i0}p_{i+1,j}+1,
\end{eqnarray*}
from Proposition~\ref{prop:unpuncturedmesh}.
Hence $d_0m_{i0}m_{i+1,0}=m_{i+1,i}+1$ as required.

\noindent {\bf Proof of (3), Case (II):}
We next assume that $i,i+1$ and $k$ are distinct while $i+1=j$.
See Figure~\ref{fig:part3caseII}.

(a) There are $|\M(\{0\}\cup V_{k,i-1},\T|_{P(i)})|$
possible restrictions to $P(i)$ of a matching in $\M_{i+1,i}$
in which $k$ is allocated a triangle in $P(i)$.
By Theorem~\ref{thm:broline} this is equal to $p_{i0}=m_{i0}$.
There are $|\M(V_{i+2,k-1},\T|_{Q(i)})|$ possible restrictions to $Q(i)$,
which by Theorem~\ref{thm:broline} is equal to
$|\M(\{0\},\T|_{Q(i)})|=d_0-1$.
We see that there are $(d_0-1)m_{i0}$ matchings in $\M_{i+1,i}$
in which $k$ is allocated a triangle in $P(i)$.

(b) There are $|\M(V_{k+1,i-1},\T|_{P(i)})|$ possible restrictions to $P(i)$
of a matching in $\M_{i+1,i}$ in which $k$ is allocated a triangle in $Q(i)$.
By Theorem~\ref{thm:broline} this is equal to $p_{ki}$.
There are $|\M(V_{i+2,k},\T|_{Q(i)})|$ possible restrictions to $Q(i)$,
which by Theorem~\ref{thm:broline} is equal to $q_{0,i+1}$. This
equals $1$ by Remark~\ref{rem:unpuncturedone}.
We see that there are $p_{ki}$ matchings in $\M_{i+1,i}$ in which $k$
is allocated a triangle in $Q(i)$.

By Lemma~\ref{lem:decomposition} we obtain
$m_{i+1,i}=(d_0-1)m_{i0}+p_{ki}.$
But by Proposition~\ref{prop:unpuncturedmesh}, we have that
$p_{i0}=p_{ki}+1$ (using the fact that $D_{i+1,k}$ is an arc in
$\T|_{P(i)}$ and Remark~\ref{rem:unpuncturedone}), so we get that
$m_{i+1,i}=d_0m_{i0}-1$, so $m_{i+1,i}+1=d_0m_{i0}m_{i+1,0}$ as required,
since $m_{i+1,0}=1$ by Lemma~\ref{lem:onecase}.

\noindent {\bf Proof of (3), Case (III):}
We next assume that $i,i+1$ and $j$ are distinct while $i=k$.
See Figure~\ref{fig:part3caseIII}.

\begin{figure}\label{fig:part3caseIIIandIV}
\begin{center}
\subfigure[Case III]
   {\includegraphics[scale=.4]{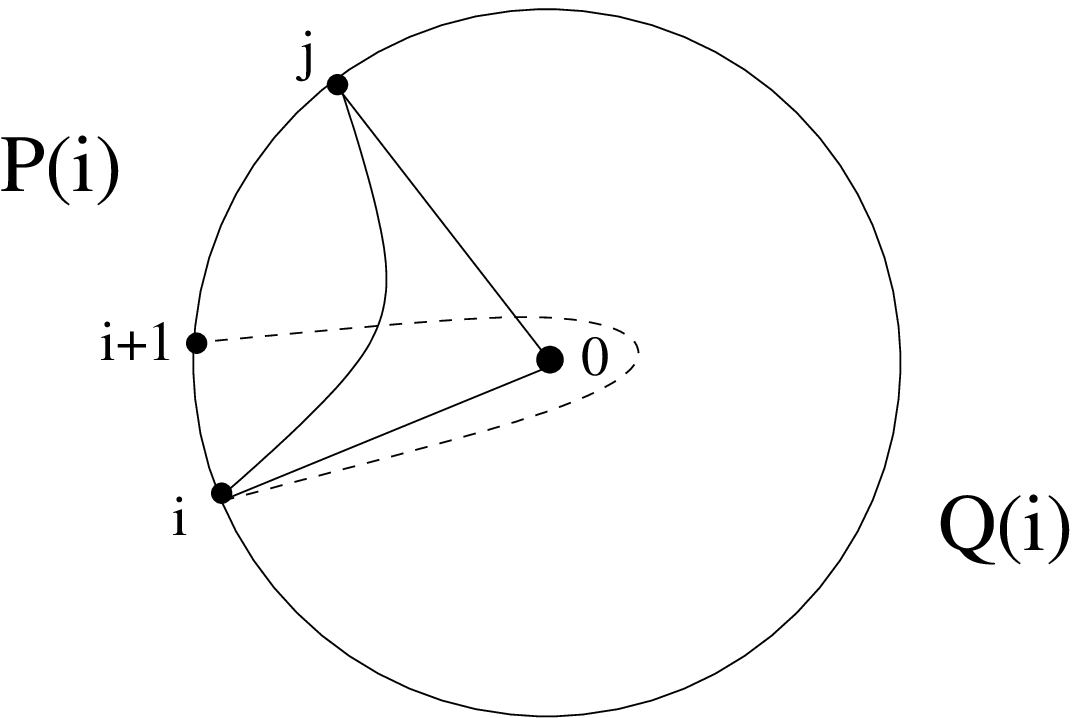}\label{fig:part3caseIII}}
\hskip 2cm
\subfigure[Case IV]
   {\includegraphics[scale=.4]{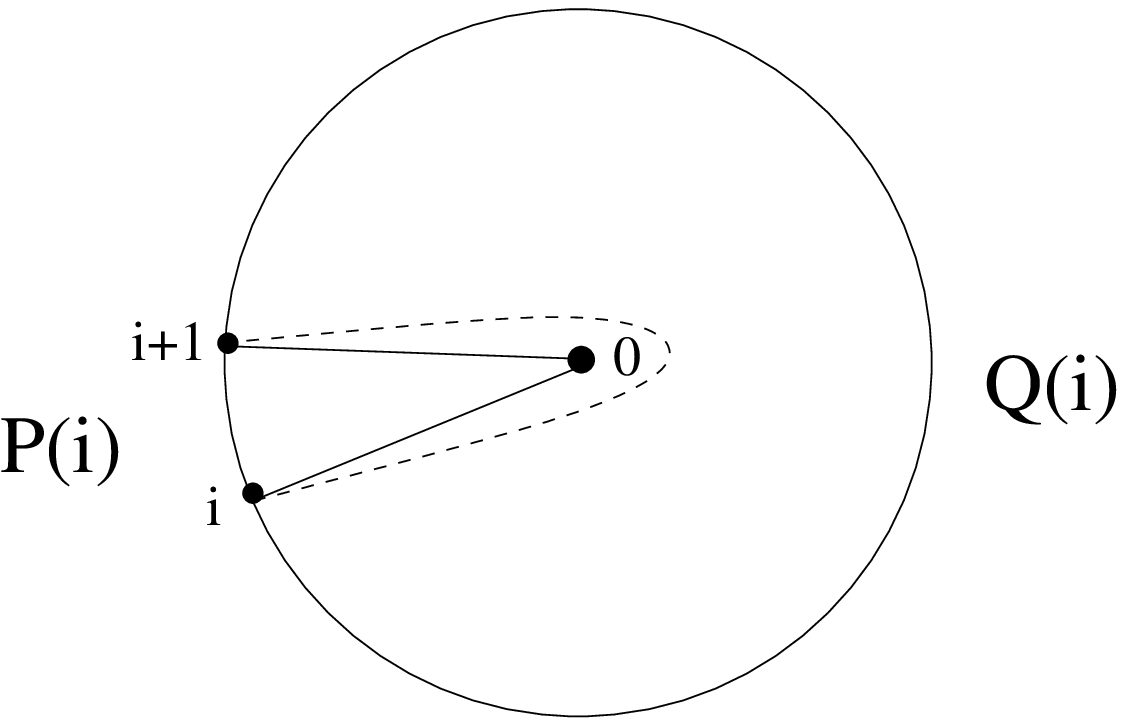}\label{fig:part3caseIV}}
\end{center}
\caption{Proof of Proposition~\ref{prop:relations}(3), Cases (III) and (IV)}
\end{figure}

It follows from symmetry with Case (II) above that
$d_0m_{i+1,0}=m_{i+1,i}+1$, so $d_0m_{i+1,0}m_{i0}=m_{i+1,i}+1$
as required, since $m_{i0}=1$ by Lemma~\ref{lem:onecase}.

\noindent {\bf Proof of (3), Case (IV):}
We next assume that $i=j$ and $i+1=k$.
See Figure~\ref{fig:part3caseIV}.

A matching in $\M_{i+1,i}$ is determined by its restriction to $Q(i)$.
The number of possible restrictions is
$|\M(V_{i+2,i-1},\T|_{Q(i)})|$
which is $|\M(\{0\},\T|_{Q(i)})|=d_0-1$ by Theorem~\ref{thm:broline}.
So $m_{i+1,i}=d_0-1$. Since $m_{i0}=m_{i+1,0}=1$ by Lemma~\ref{lem:onecase}
we obtain $m_{i+1,i}+1=d_0m_{i0}m_{i+1,0}$ as required.

\noindent {\bf Proof of (3), Case (V):}
We next assume that $i,i+1$ and $j$ are distinct, with $j=k$.
See Figure~\ref{fig:part3caseV}.
\begin{figure}
\begin{center}
\includegraphics[scale=.4]{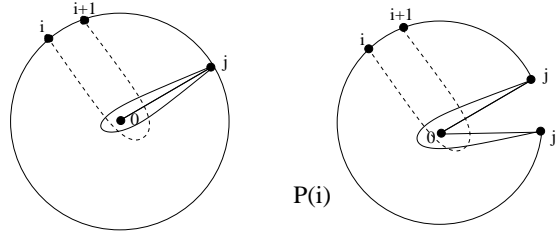}
\caption{Proof of Proposition~\ref{prop:relations}(3), Case (V).}
\label{fig:part3caseV}
\end{center}
\end{figure}

There are $|\M(V_{i+2,j},\T|_{P(i)})|\cdot |\M(V_{j+1,i-1},\T|_{P(i)})|$
possible restrictions to $P(i)$ of a matching in $\M_{i+1,i}$ in which $j$
gets allocated a triangle, since $D_{i+1,i}$ splits $P(i)$ completely.
There are
$|\M(V_{i+2,j-1},\T|_{P(i)})|\cdot |\M(V_{j',i-1},\T|_{P(i)})|$
possible restrictions to $P(i)$ in which $j'$ gets allocated a triangle.
Thus the total number of matchings in $\M_{i+1,i}$ is the sum of these,
which by Theorem~\ref{thm:broline} is equal to
$p_{i+1,0}p_{j'i}+p_{i+1,j}p_{0i}$.
By Proposition~\ref{prop:unpuncturedmesh}, $p_{i0}p_{jj'}=p_{ij}p_{0j'}+p_{ij'}p_{0j}$,
so $p_{i0}=p_{ij}+p_{ij'}$. We thus have that
\begin{eqnarray*}
m_{i+1,i} &=& m_{i+1,0}(p_{i0}-p_{ij})+m_{i0}p_{i+1,j} \\
&=& m_{i+1,0}m_{i0}-m_{i+1,0}p_{ij}+m_{i0}p_{i+1,j} \\
&=& m_{i+1,0}m_{i0}-1 \\
&=& d_0m_{i0}m_{i+1,0}-1,
\end{eqnarray*}
as required. Here we use the fact that
$$p_{i+1,0}p_{ij}=p_{i+1,j}p_{i0}+p_{0j}p_{i,i+1}=p_{i0}p_{i+1,j}+1$$
from Theorem~\ref{thm:broline} and Remark~\ref{rem:unpuncturedone},
and the fact that $d_0=1$.

\noindent {\bf Proof of (3), Case (VI):}
We finally assume that $i=j=k$.
See Figure~\ref{fig:part3caseVI}.
\begin{figure}
\begin{center}
\includegraphics[scale=.4]{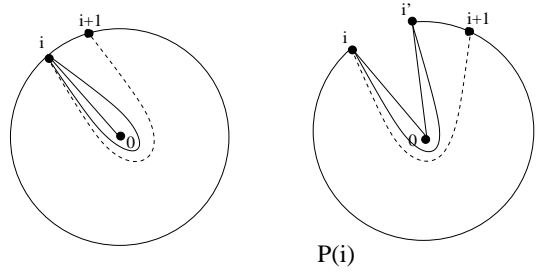}
\caption{Proof of Proposition~\ref{prop:relations}(3), Case (VI).}
\label{fig:part3caseVI}
\end{center}
\end{figure}

Since a matching in $\M_{i+1,i}$ is determined its restriction to
$P(i)$, which is a matching between $\T|_{P(i)}$ and $V_{i+2,\ldots ,i-1}$
we have that $m_{i+1,i}=|\M(V_{i+2,i-1},\T|_{P(i)})$
which equals $p_{i+1,i}$ by Theorem~\ref{thm:broline}.
By Proposition~\ref{prop:unpuncturedmesh}, we have that
$$p_{i+1,0}p_{i'i}=p_{i+1,i}p_{0i'}+p_{i',i+1}p_{i0},$$
so, since $m_{i+1,0}=p_{i+1,0}$ (by definition of $m_{i+1,0}$)
and $p_{i'i}=p_{0i'}=p_{i0}=1$ (using Remark~\ref{rem:unpuncturedone}),
we obtain $m_{i+1,0}=m_{i+1,i}+1$, so $d_0m_{i0}m_{i+1,0}=m_{i+1,i}+1$,
as required, since $d_0=m_{i0}=1$.

\noindent {\bf Proof of (4)}: We note that $m_{i0}m_{i+1,i+1}=
m_{i0}d_0m_{i+1,0}=m_{ii}m_{i+1,0}$ by Lemma~\ref{lem:degeneratecase},
so (4) follows from (3).

The proposition is proved.
\end{proof}

\begin{theorem} \label{thm:dnfriezepattern}
Let $\T$ be a triangulation of $S_{\odot}$.
Then the numbers $m_{ij}$, $m_{ii}$ and $m_{i0}$ in
Definition~\ref{def:matchingnumbers} (arranged as in the introduction)
form a frieze pattern of type $D_N$.
\end{theorem}

\begin{proof}
This follows immediately from Proposition~\ref{prop:relations} and
Lemma~\ref{lem:positivenumbers}.
\end{proof}

%
\section{Tagged triangulations}
\label{sec:tagged}
%

In Section~\ref{sec:mesh-rel}, we have established a way to 
obtain a frieze pattern of type $D_N$ from a triangulation 
of a punctured disc with $N$ boundary vertices 
using the matching numbers. Here, we will show that there 
exist frieze patterns (of type $D_N$) which cannot be obtained 
in the same way. We will describe another way to construct 
such frieze patterns and will show how they can be associated 
to tagged triangulations (as defined below).

\begin{definition}\label{def:pattern}
We denote the frieze pattern of type $D_N$ associated 
in Theorem~\ref{thm:dnfriezepattern}
to the triangulation $\T$ of a punctured disc by $F(\T)$. 
\end{definition}

The frieze patterns of the form $F(\T)$ do not give 
all possible frieze patterns of type $D_N$, as we will see now. 

\begin{definition}\label{def:switching}
Let $P$ be a frieze pattern of type $D_N$. We define $\iota(P)$ 
to be the pattern obtained from $P$ by interchanging the last 
two rows. 
\end{definition}

Thus, if $\T$ is a triangulation of $S_{\odot}$, $\iota(F(\T))$ is 
also a frieze pattern of type $D_N$. 
If a triangulation has only one triangle 
at the puncture (i.e. $d_0=1$), then $m_{ii}=m_{i0}$ for all 
$i$ (by Lemma~\ref{lem:degeneratecase}), 
so the process of interchanging the last two rows does not 
give a new pattern, so $\iota (F(\T))=F(\T)$. 
We claim that in general, $\iota (F(\T))$ cannot be obtained 
via a triangulation of a punctured disc: 

\begin{remark}
Let $\T$ be a triangulation of a punctured disc $S_{\odot}$. 
If $\T$ has at least two central arcs then there is 
no triangulation $\T'$ of $S_{\odot}$ such that $\iota (F(\T))=F(\T')$. 
\end{remark}

\begin{proof}
Let $F(\T)$ be obtained from the matching numbers $m_{ij}$ 
of $\T$ and let $\iota (F(\T))$ be defined as above. 
Assume that $d_0>1$. Then since $m_{i0}=d_0m_{ii}$ 
(Lemma~\ref{lem:degeneratecase}) we have 
$m_{i0}<m_{ii}$. Now 
if there exists a triangulation $\T'$ of $S_{\odot}$ 
with matching numbers 
$m'_{ij}$ giving rise to $\iota (F(\T))$, then by 
Lemma~\ref{lem:degeneratecase} (applied to 
$\T'$), the matching numbers of $\T'$ must 
satisfy $m'_{i0}<m'_{ii}$ for all $i$. 
But by construction, $m'_{i0}=m_{ii}>m_{i0}=m_{ii}'$. 
\end{proof} 

\begin{figure}
\begin{center}
\includegraphics[scale=.4]{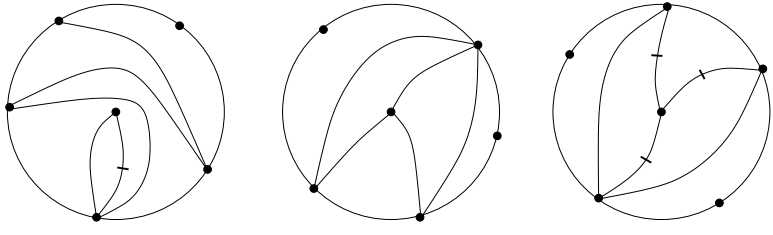}
\caption{Three tagged triangulations}
\label{fig:tagged}
\end{center}
\end{figure}

We thus have a way to construct a second frieze pattern for 
every triangulation with $d_0>1$. 
We now recall the definition of tagged arcs resp. of tagged 
edges as introduced in~\cite[Definition 7.1]{fst06} for triangulated 
surfaces and in ~\cite[Section 2.1]{schiffler06} for punctured polygons 
respectively. The idea is to attach labels to arcs in the punctured 
disc $S_{\odot}$: 

\begin{definition} \label{def:taggedarcs}
Let $S_{\odot}$ be a punctured disc with $N$ boundary vertices 
with puncture $0$. The set of {\em tagged arcs of $S_{\odot}$} 
is the following: 
\[
\{ D_{ij}^1,\ D_{i0}^1,\ D_{i0}^{-1} \mid 1\le i,j\le N,\ j\neq i,i+1 \}.
\]
\end{definition}
We will often just write $D_{ij}$ instead of $D_{ij}^1$. 
When drawing tagged arcs, we indicate the label $-1$ of an arc by a short crossing 
line on it. Arcs labelled $1$ are drawn as the usual arcs, cf. Figure~\ref{fig:tagged}. 

\begin{remark}\label{rem:arcbijection}
Note that the tagged arcs are in bijection with the arcs in the sense of 
Definition~\ref{def:arc}: 
The arcs $D_{ij}^1$ ($j\neq i,i+1$) correspond to the usual arcs $D_{ij}$ ($j\neq i,i+1$), 
the central arcs $D_{i0}^1$ to the usual central arcs $D_{i0}$ and the 
arcs $D_{i0}^{-1}$ of label $-1$ correspond to the loops $D_{ii}$. 
\end{remark}

Following~\cite[Definition 7.7]{fst06} and~\cite[Section 2.4]{schiffler06} we can now define tagged 
triangulations of punctured discs. 

\begin{definition}\label{def:taggedtriangulation}
A {\em tagged triangulation} $\widetilde{\T}$ of $S_{\odot}$ is a collection of 
tagged arcs of $S_{\odot}$ obtained from a triangulation $\T$ of $S_{\odot}$ 
as follows: 

\begin{tabular}{ll}
(i) & If $d_0=1$, replace the pair $D_{ii}$, $D_{i0}$ by the tagged arcs $D_{i0}^{-1}$, $D_{i0}^1$\\ 
    & and replace the arcs $D_{ij}$ by tagged arcs $D_{ij}^1$. \\
(ii) & If $d_0=k\ge 2$, replace all central arcs $D_{i_j,0}$ either by $k$ tagged arcs $D_{i_j,0}^1$ \\
     &  or by $k$ tagged arcs $D_{i_j,0}^{-1}$. 
\end{tabular}
\end{definition}

See Figure~\ref{fig:tagged} for three 
tagged triangulations for a punctured disc with $5$ boundary vertices. 

In other words: if $\T$ is a triangulation of $S_{\odot}$ containing a unique central 
arc $D_{i0}$ and hence also the loop $D_{ii}$ (in particular, $d_0=1$), then $\T$ gives rise to a 
tagged triangulation 
$\widetilde{\T}$ whose only central tagged arcs are $D_{i0}^1$ and $D_{i0}^{-1}$. 
If $\T$ has central arcs $D_{i_1,0},\dots,D_{i_k,0}$, $k\ge 2$ (so $d_0=k\ge 2$), 
then 
$\T$ gives rise to two tagged arc triangulations: one with $k$ central tagged 
arcs $D_{i_j,0}^1$ and one with $k$ central tagged arcs $D_{i_j,0}^{-1}$. 

\begin{definition} \label{def:Ftagged}
Let $\widetilde{\T}$ be an arbitrary tagged triangulation of a
punctured disc $S_{\odot}$. 
Then we associate a frieze pattern $F(\widetilde{\T})$ to it
as follows:

(a) If $\widetilde{\T}$ has only arcs labelled $1$, let $\T$ be the triangulation
of $S_{\odot}$ obtained via the bijection of Remark~\ref{rem:arcbijection}. We
then define $F(\widetilde{\T})$ as the frieze pattern of the
matching numbers of $\T$, i.e. we set $F(\widetilde{\T})=F(\T)$.

(b) If $\widetilde{\T}$ has exactly one tagged arc labelled $-1$, say $D_{i0}^{-1}$ and
hence also $D_{i0}^1$, then let $\T$ be the triangulation of $S_{\odot}$
obtained through the bijection in Remark~\ref{rem:arcbijection}. As in (a),
we let $F(\widetilde{\T})$ be the frieze pattern $F(\T)$ containing the
matching numbers of $\T$.

(c) Suppose that $\widetilde{\T}$ has at least two arcs labelled $-1$.
Then every arc incident with $0$ must be labelled $-1$.
Let $\T'$ be the triangulation of $S_{\odot}$ consisting of all
arcs $D_{i0}$ such that $D_{i0}^{-1}$ lies in $\widetilde{\T}$ and all
arcs $D_{ij}$ (for $i$ and $j$ boundary vertices of $S_{\odot}$) such that
$D_{ij}^1$ lies in $\widetilde{\T}$.
Then we set $F(\widetilde{\T})=\iota (F(\T'))$.
\end{definition}

We also note the following:

\begin{prop}
The entries in a frieze pattern $M$ of type $D_N$ associated to a
tagged triangulation $\widetilde{\T}$ of a punctured disc are determined by
the numbers of triangles incident with the marked points, together with the
number of triangles incident with the puncture.
\end{prop}

\begin{proof}
Suppose first that $\T$ is a triangulation (without tags).
We remark that it follows from the definition of $F(\T)$ that
an entry $m_{i,i+2}$ in the first row (with $1\leq i\leq N$) is equal to
the number of triangles incident with vertex $i$ in $\T$. 

It is clear that,
by induction on the rows, these entries determine the entries
$m_{ij}$ in the pattern for all boundary vertices $i,j$ with $i\not=j$, via
relation~\eqref{eqn:relation1prop} in Proposition~\ref{prop:relations}.
Relations~\eqref{eqn:relation2prop} and~\eqref{eqn:relation3prop} can be used
to determine $m_{ii}m_{i0}$ for each boundary vertex $i$.
Using Lemma~\ref{lem:degeneratecase} we see that all entries in $F(\T)$
are determined in cases (a) and (b) of Definition~\ref{def:Ftagged}.
In case (c) we have $F(\widetilde{\T})=\iota (F(\T))$ 
and it is clear that the same approach works for such $\widetilde{\T}$. 
\end{proof}

%
\section{Cluster algebras and frieze patterns of type $D_N$}
\label{sec:clusteralgebras}
%

Let $\mathcal{A}$ be a cluster algebra~\cite{fominzelevinsky02} of type
$D_N$, as in the classification of cluster algebras of finite
type~\cite{fominzelevinsky03}. We consider the case in which all coefficients
are set to $1$. The algebra $\mathcal{A}$ is a subring of the rational
function field $\mathbb{F}:=\mathbb{Q}(x_1,x_2,\ldots ,x_N)$.
It is determined by an initial \emph{seed},
i.e. a pair consisting of a transcendence basis
$\mathbf{x}_0$ of $\mathbb{F}$ over $\mathbb{Q}$ and a
skew-symmetric integer matrix $B_0$ with rows and columns
indexed by $\mathbf{x}_0$. For each element of $\mathbf{x}_0$,
a new seed can be produced from $(\mathbf{x}_0,B_0)$ by a combinatorial
process known as \emph{mutation}. We obtain a collection of seeds via
arbitrary iterative mutation of $(\mathbf{x}_0,B_0)$.

The transcendence bases arising are known as a \emph{clusters} and
$\mathcal{A}$ is generated by their union. The generators are known as
\emph{cluster variables}. We note that a skew-symmetric matrix $B$
appearing in a seed $(\mathbf{x},B)$ can be encoded by a quiver,
with vertices indexed by $\mathbf{x}$ and $b_{xy}$ arrows from the
vertex indexed by $x$ to the vertex indexed by $y$ whenever $b_{xy}>0$.

By~\cite[7.11]{fst06} (see also~\cite{schiffler06}), the cluster variables of
$\mathcal{A}$ are in bijection with the set of all tagged arcs in $S_{\odot}$,
and the seeds of $\mathcal{A}$ are in bijection with the tagged
triangulations of $S_{\odot}$.
The matrix $B$ of a seed corresponding to a given tagged triangulation
is described in~\cite[9.6]{fst06}.

Combining the above bijection between cluster variables and tagged arcs
with the bijection in Remark~\ref{rem:arcbijection}, we obtain a
bijection between cluster variables and the arcs of $S_{\odot}$ in the sense of
Definition~\ref{def:arc}. For $D_{ij}$ (respectively, $D_{i0}$) an arc of
$S_{\odot}$, let $x_{ij}$ (respectively, $x_{i0}$) denote the corresponding
cluster variable.

\begin{definition}
Let $\widetilde{\T}$ be a tagged triangulation of $S_{\odot}$, and let
$\mathbf{x}$ be
the corresponding cluster of $\mathcal{A}$. We know
from (a special case of)~\cite[3.1]{fominzelevinsky01} that each cluster
variable of $\mathcal{A}$ can be expressed as a Laurent polynomial in the
elements of $\mathbf{x}$ with integer coefficients.
For any arc $D_{ij}$ (respectively, $D_{i0}$) of $S_{\odot}$, let
$u_{ij}$ (respectively, $u_{i0}$) denote the integer obtained
from $x_{ij}$ (respectively, $u_{i0}$) when the elements of $\mathbf{x}$
are all specialised to $1$.

Let $F_c(\widetilde{\T})$ be the array of integers $u_{ij}$ written in
the same positions as the $m_{ij}$ in Figures~\ref{fig:dpatterneven}
and~\ref{fig:dpatternodd}.
\end{definition}

\begin{prop} \label{prop:clusterrelations}
For any tagged triangulation $\widetilde{\T}$, the array $F_c(\widetilde{\T})$
defined above is a frieze pattern of type $D_N$.
\end{prop}

\begin{proof}
If $D_{ij}$ (respectively, $D_{i0}$) lies in $\T$, then $u_{ij}=1$
(respectively, $u_{i0}=1$) is positive.
The positivity of any entry in $F_c(\T)$ then follows by
induction on the number of mutations needed from $\mathbf{x}$ to obtain the
corresponding cluster variable, using the exchange relations in $\mathcal{A}$.
The frieze relations follow from some of the exchange relations
in $\mathcal{A}$, i.e.\ those arising
from~\cite[5.1, 5.2]{schiffler06} (cf.~\cite{bmr04}).
Relation~\eqref{eqn:relation1} (see the introduction) arises from case (1)
of~\cite[5.1]{schiffler06} in the case where $a\neq d$ in Schiffler's
notation. Relation~\eqref{eqn:relation2} arises from case (1),
in the case where $a=d$.
Relations~\eqref{eqn:relation3} and~\eqref{eqn:relation4}
arise from case (2) of~\cite[5.1]{schiffler06}. (Note that the
exchange relations for a cluster algebra of type $D_N$ are also described
in~\cite[12.4]{fominzelevinsky03}).
\end{proof}

\begin{definition} \label{def:slice}
A \emph{slice} of a frieze pattern of type $D_N$ is defined as follows.
We initially select an entry $E_1$ in the top row.
For $2\leq i\leq N-1$, suppose that an entry $E_{i-1}$ in the
$i-1$st row has already been chosen. We select an entry $E_i$ in the
$i$th row which is either immediately down and to the right of $E_{i-1}$
or immediately down and to the left of $E_{i-1}$.
In addition, we select an entry $E_N$ in row $N$ below the entry
immediately to the right of and below $E_{N-2}$ or below the entry\
immediately to the left of and below $E_{N-2}$. 
\end{definition}

Definition~\ref{def:slice} is motivated by the notion of a slice in
representation theory (see~\cite{ringel84}).

\begin{remark}
Let $(\mathbf{x},B)$ be a seed of $\mathcal{A}$. Let $\Gamma$ be the
quiver associated to $B$ and let $\widetilde{\T}$ be a tagged triangulation.
Then it is easy to check using~\cite[9.6]{fst06} that
$\Gamma$ is an orientation of the Dynkin diagram of type $D_N$ if and only
if the corresponding subset of $F_c(\widetilde{\T})$ is a slice in
the above sense.
\end{remark}

\begin{lemma} \label{lem:slice}
The entries in a slice of a frieze pattern of type $D_N$ determine
the entire pattern.
\end{lemma}

\begin{proof}
It is straightforward to see that the frieze relations can be used to
determine all of the entries in the pattern.
\end{proof}

\begin{theorem} \label{thm:clustertypeD}
Suppose that the cluster $\mathbf{x}$ corresponds to a slice in the
frieze pattern, and let $\widetilde{\T}$ be the corresponding tagged
triangulation. Then the frieze patterns $F(\widetilde{\T})$ and
$F_c(\widetilde{\T})$ coincide.
\end{theorem}

\begin{proof}
As stated in Definition~\ref{def:Ftagged}, $F(\widetilde{\T})$ is
a frieze pattern of type $D_N$. By Proposition~\ref{prop:clusterrelations},
$F_c(\widetilde{\T})$ is a frieze pattern of type $D_N$.
Let $x_{ij}$ (respectively, $x_{i0}$) be the cluster variable corresponding to
a tagged arc of $\widetilde{\T}$. Then $u_{ij}=1$ (respectively, $u_{i0}=1$) by
definition. If $\widetilde{\T}$ has at most one arc labelled $-1$ then
$m_{ij}=1$ (respectively $m_{i0}=1$) by Lemma~\ref{lem:onecase}.
If $\widetilde{\T}$ has more than one arc tagged with $-1$ then $m_{ij}=1$
(respectively, $m_{i0}=1$) by
the definition of $F(\widetilde{\T})$
((c) in Definition~\ref{def:Ftagged}).

The fact that $F(\widetilde{\T})$ and $F_c(\widetilde{\T})$ coincide then
follows from iterative application of the relations of
Proposition~\ref{prop:relations} and the frieze relations for
$F_c(\widetilde{\T})$ (using Proposition~\ref{prop:clusterrelations}).
We use Lemma~\ref{lem:slice}.
\end{proof}

Finally, motivated by the situation for classical frieze patterns,
we make the following conjectures:

\begin{conjecture}
Let $\widetilde{\T}$ be any tagged triangulation of $S_{\odot}$.
Then $F(\widetilde{\T})$ and $F_c(\widetilde{\T})$ coincide.
\end{conjecture}

\begin{conjecture}
Every frieze pattern of type $D_N$ is of the form $F(\widetilde{\T})$ for
some tagged triangulation $\widetilde{\T}$ of $S_{\odot}$.
\end{conjecture}

%
\section{An Example}
\label{sec:example}
%

\begin{figure}
\begin{center}
\includegraphics[scale=.4]{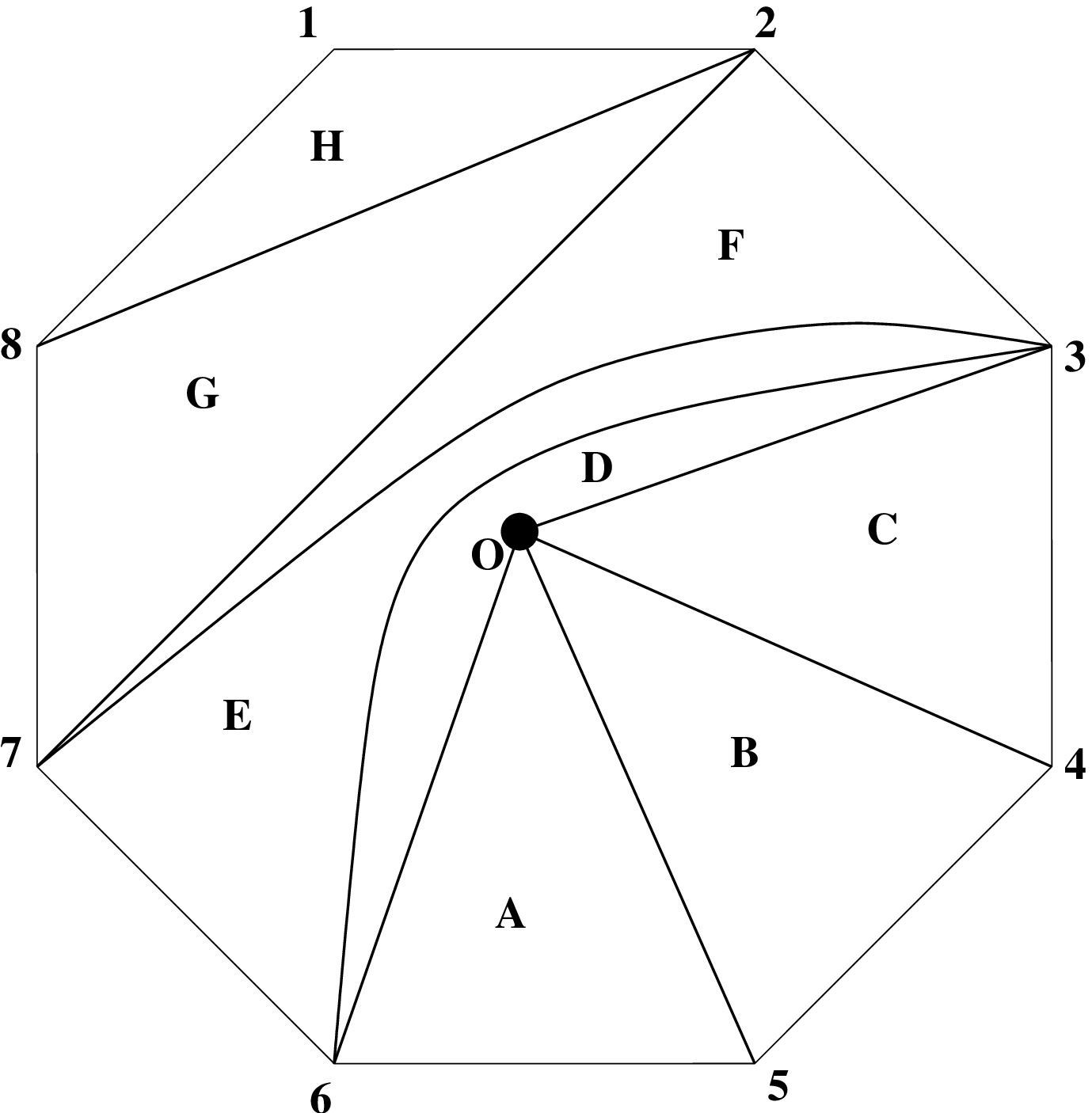}
\caption{A triangulation of a punctured octagon}
\label{D8exampletri}
\end{center}
\end{figure}

We consider an example of type $D_8$. We take the triangulation $\T$ of
the punctured disc with $8$ marked points on its boundary displayed in
displayed in Figure~\ref{D8exampletri}. Each triangle has been
labelled with a letter for reference. The corresponding frieze
pattern is shown in Figure~\ref{D8examplefrieze}.
Each entry in the frieze pattern (apart from the first row) is labelled by
its corresponding diagonal.

For example, the entry under $D_{27}$ is $23$, since there are $23$ matchings
between $\T|_{S_{27}}$ and the vertices $V_{3,6}=\{3,4,5,6\}$.
Note that the triangles $D$ and $E$ are split by the arc $D_{27}$ and
we include in the count matchings in which both resulting regions are
allocated to vertices. Thus, for example, the matching
$3E\ 4B\ 5A\ 6E$ is included.

The entry under $D_{52}$ is $5$, since there are $5$ matchings
between $\T|_{S_{52}}$ and the vertices $V_{5,2}=\{6,7,8,1\}$.
Note that in this case the arc $S_{52}$ does not split any triangles.
The matchings are:
$6A\ 7E\ 8G\ 1H$, $6D\ 7E\ 8G\ 1H$, $6A\ 7F\ 8G\ 1H$, $6D\ 7F\ 8G\ 1H$,
and $6E\ 7F\ 8G\ 1H$.

The entry under $D_{22}$ is $12$, since there are $12$ matchings between
$\T$ and the vertices $V_{31}=\{3,4,5,6,7,8,1\}$ and $0$.
In this case each triangle can only be allocated to a single vertex, even
triangle $F$ which is split by the arc $D_{22}$.
The $12$ matchings are shown in Figure~\ref{D8examplematchings}.

Finally, we note that for all $i$, $D_{ii}/D_{i0}=4$, the number of
triangles of $\T$ incident with $0$.

\begin{figure}
\begin{center}
\setlength{\arraycolsep}{0mm}
$$
\begin{array}{cccccccccccccccccccccccc}
1 && 1 && 1 && 1 && 1 && 1 && 1 && 1 && 1 &&&& \\
&D_{13}&      &D_{24}&      &D_{35}&      &D_{46}  &      &D_{57}  &      &D_{68}  &      &D_{71}  &      &D_{82}  &      &D_{13}  &      &        &      &        &      & \\
&3     &      &4     &      &2     &      &2       &      &3       &      &3       &      &2       &      &1       &      &3       &      &        &      &        &      & \\
&      &D_{14}&      &D_{25}&      &D_{36}&        &D_{47}&        &D_{58}&        &D_{61}&        &D_{72}&        &D_{83}&        &D_{14}&        &      &        &      & \\
&      &11    &      &7     &      &3     &        &5     &        &8     &        &5     &        &1     &        &2     &        &11    &        &      &        &      & \\
&      &      &D_{15}&      &D_{26}&      &D_{37}  &      &D_{48}  &      &D_{51}  &      &D_{62}  &      &D_{73}  &      &D_{84}  &      &D_{15}  &      &        &      & \\
&      &      &19    &      &10    &      &7       &      &13      &      &13      &      &2       &      &1       &      &7       &      &19      &      &        &      & \\
&      &      &      &D_{16}&      &D_{27}&        &D_{38}&        &D_{41}&        &D_{52}&        &D_{63}&        &D_{74}&        &D_{85}&        &D_{16}&        &      & \\
&      &      &      &27    &      &23    &        &18    &        &21    &        &5     &        &1     &        &3     &        &12    &        &27    &        &      & \\
&      &      &      &      &D_{17}&      &D_{28}  &      &D_{31}  &      &D_{42}  &      &D_{53}  &      &D_{64}  &      &D_{75}  &      &D_{86}  &      &D_{17}  &      & \\
&      &      &      &      &62    &      &59      &      &29      &      &8       &      &2       &      &2       &      &5       &      &17      &      &62      &      & \\
&      &      &      &      &      &D_{18}&        &D_{21}&        &D_{32}&        &D_{43}&        &D_{54}&        &D_{65}&        &D_{76}&        &D_{87}&        &D_{18}& \\
&      &      &      &      &      &159   &        &95    &        &11    &        &3     &        &3     &        &3     &        &7     &        &39    &        &159   & \\
&      &      &      &      &      &      &D_{11}&      &D_{20}&      &D_{33}&      &D_{40}&      &D_{55}&      &D_{60}&      &D_{77}&      &D_{80}&      &D_{11} \\
&      &      &      &      &      &      &32      &      &3       &      &4       &      &1       &      &4       &      &1       &      &8       &      &5       &      &32       \\
&      &      &      &      &      &      &D_{10}&      &D_{22}&      &D_{30}&      &D_{44}&      &D_{50}&      &D_{66}&      &D_{70}&      &D_{88}&      &D_{10} \\
&      &      &      &      &      &      &8       &      &12      &      &1       &      &4       &      &1       &      &4       &      &2       &      &20      &      &8         
\end{array}
$$
\end{center}
\caption{A frieze pattern of type $D_8$}
\label{D8examplefrieze}
\end{figure}

\begin{figure}
\begin{center}
\noindent
$3C\ 4B\ 5A\ 6E\ 7F\ 8G\ 1H\ OD$, \\
$3D\ 4B\ 5A\ 6E\ 7F\ 8G\ 1H\ OC$, \\
$3D\ 4C\ 5A\ 6E\ 7F\ 8G\ 1H\ OB$, \\
$3D\ 4C\ 5B\ 6E\ 7F\ 8G\ 1H\ OA$, \\
$3E\ 4B\ 5A\ 6D\ 7F\ 8G\ 1H\ OC$, \\
$3E\ 4C\ 5A\ 6D\ 7F\ 8G\ 1H\ OB$, \\
$3E\ 4C\ 5B\ 6A\ 7F\ 8G\ 1H\ OD$, \\
$3E\ 4C\ 5B\ 6D\ 7F\ 8G\ 1H\ OA$, \\
$3F\ 4B\ 5A\ 6D\ 7E\ 8G\ 1H\ OC$, \\
$3F\ 4C\ 5A\ 6D\ 7E\ 8G\ 1H\ OB$, \\
$3F\ 4C\ 5B\ 6A\ 7E\ 8G\ 1H\ OD$, \\
$3F\ 4C\ 5B\ 6D\ 7E\ 8G\ 1H\ OA$.
\end{center}
\caption{Matchings for the arc $D_{22}$}
\label{D8examplematchings}
\end{figure}


\begin{thebibliography}{CoCo73b}

\bibitem[BCI74]{bci74}
D. Broline, D. W. Crowe and I. M. Isaacs.
\textit{The geometry of frieze patterns.}
Geom. Ded. \textbf{3} (1974), 171--176.

\bibitem[BMR04]{bmr04}
A. B. Buan, R. J. Marsh and I. Reiten.
\textit{Cluster mutation via quiver representations.}
Preprint, arXiv:math.RT/0412077, 2004, to appear in Comm. Math. Helv.

\bibitem[B99]{burman99}
Y.M. Burman.
\textit{Triangulations of disc with boundary and the homotopy principle for
functions without critical points.}
Annals of Global Analysis Geometry,
\textbf{17} (1999), 221--238.

\bibitem[CaCh06]{calderochapoton06}
P. Caldero and F. Chapoton.
\textit{Cluster algebras as Hall algebras of quiver representations.}
Commentarii Mathematici Helvetici, 81, (2006), 595-616.

\bibitem[CaPr03]{carrollprice}
G. Carroll and G. Price.
\textit{Two new combinatorial models for the Ptolemy recurrence.}
Unpublished memo (2003).

\bibitem[CoCo73a]{cc73a}
J. H. Conway and H. S. M. Coxeter,
\textit{Triangulated discs and frieze patterns.}
Math. Gaz. \textbf{57} (1973), 87--94.

\bibitem[CoCo73b]{cc73b}
J. H. Conway and H. S. M. Coxeter,
\textit{Triangulated discs and frieze patterns.}
Math. Gaz. \textbf{57} (1973), 175--186.

\bibitem[FST06]{fst06}
S. Fomin, D. Shapiro and D. Thurston,
\textit{Cluster algebras and triangulated discs. Part I: Cluster complexes.}
Preprint arxiv:math.RA/0608367v3, August 2006.

\bibitem[FZ01]{fominzelevinsky01}
S. Fomin and A. Zelevinsky.
\textit{$Y$-systems and generalized associahedra.}
Ann. Math. \textbf{158} (2003), no. 3,  977--1018.

\bibitem[FZ02]{fominzelevinsky02}
S. Fomin and A. Zelevinsky.
\textit{Cluster algebras I: Foundations.}
J. Amer. Math. Soc. 15 (2002), no. 2, 497--529.

\bibitem[FZ03]{fominzelevinsky03}
S. Fomin and A. Zelevinsky,
\textit{Cluster algebras II: Finite type classification.}
Invent. Math. \textbf{154} (2003), no. 1, 63--121.

\bibitem[G85]{grimaldi85}
R. P. Grimaldi,
\textit{Discrete and Combinatorial Mathematics: An Applied Introduction.}
Addison-Wesley, 1985.

\bibitem[K04]{kuo04}
E. Kuo,
\textit{Applications of graphical condensation for enumerating matchings
and tilings.}
Theoret. Comput. Sci. \textbf{319} (2004), 29--57.

\bibitem[M07]{musiker}
G. Musiker, 
\textit{A graph theoretic expansion formula for cluster algebras of type
$B_n$ and $D_n$.}
Preprint arXiv:0710.3574v1 [math.CO], 2007.

\bibitem[P05]{propp}
J. Propp.
\textit{The combinatorics of frieze patterns and Markoff numbers.}
Preprint arxiv:math.CO/0511633, 2005.

\bibitem[R84]{ringel84}
C. M. Ringel.
\textit{Tame algebras and integral quadratic forms.}
Lecture Notes in Mathematics \textbf{1099}, Springer, Berlin, 1984.

\bibitem[S06]{schiffler06}
R. Schiffler,
\textit{A geometric model for cluster categories of type $D_n$.}
Preprint arxiv:math.RT/0608264, 2006.

\end{thebibliography}
\end{document}